\documentclass[11pt]{amsart}
\usepackage{amssymb,latexsym}
\newdimen\AAdi%
\newbox\AAbo%
\def\AAk#1#2{\s_etbox\AAbo=\hbox{#2}\AAdi=\wd\AAbo\kern#1\AAdi{}}%
\def\AAr#1#2#3{\s_etbox\AAbo=\hbox{#2}\AAdi=\ht\AAbo\raise#1\AAdi\hbox{#3}}%
\font\tenmsb=msbm10 at 12pt
\font\sevenmsb=msbm7 at 8pt
\font\fivemsb=msbm5 at 6pt
\newfam\msbfam
\textfont\msbfam=\tenmsb
\scriptfont\msbfam=\sevenmsb
\scriptscriptfont\msbfam=\fivemsb
\def\Bbb#1{{\tenmsb\fam\msbfam#1}}
\begin{document}
\newtheorem{thm}{Theorem}
\newtheorem{lem}{Lemma}
\newtheorem{cor}{Corollary}
\newtheorem{rem}{Remark}
\newtheorem{pro}{Proposition}
\newtheorem{defi}{Definition}

\newcommand{\noi}{\noindent}
\newcommand{\dis}{\displaystyle}
\newcommand{\mint}{-\!\!\!\!\!\!\int}
\newcommand{\ba}{\begin{array}}
\newcommand{\ea}{\end{array}}
\def\bx{\hspace{2.5mm}\rule{2.5mm}{2.5mm}} \def \vs{\vspace*{0.2cm}}
\def\hs{\hspace*{0.6cm}}
\def\ds{\displaystyle}
\def\p{\partial}
\def\O{\Omega}
\def\o{\omega}
\def\b{\beta}
\def\m{\mu}
\def\l{\lambda}
\def\L{\Lambda}
\def\ul{u_\lambda}
\def\D{\Delta}
\def\d{\delta}
\def\s{\sigma}
\def\e{\varepsilon}
\def\a{\alpha}
\def\tf{\tilde{f}}
\def\cqfd{%
\mbox{ }%
\nolinebreak%
\hfill%
\rule{2mm} {2mm}%
\medbreak%
\par%
}

\def\pr {\noindent {\it Proof.} }
\def\rmk {\noindent {\it Remark.} }
\def\esp {\hspace{4mm}}
\def \dsp {\hspace{2mm}}
\def \ssp {\hspace{1mm}}

\def\energy{\int_{\R^n}\u }
\def\sk{\s_k}
\def\mo{\mu_k}
\def\cal{\mathcal}
\def\I{{\cal I}}
\def\J{{\cal J}}
\def\K{{\cal K}}
\def \OM{\overline{M}}
\def\fk{{{\cal F}}_k}
\def\M1{{{\cal M}}_1}
\def\Fk{{\cal F}_k}
\def\Fl{{\cal F}_l}
\def\En{{\cal E}_{n/2}}
\def\F{\cal F}
\def\FF{\cal F}
\def\tFF{\tilde {\cal F}}
\def\Gk{{\Gamma_k^+}}
\def\n{\nabla}
\def\uuu{{\n ^2 u+du\otimes du-\frac {|\n u|^2} 2 g_0+S_{g_0}}}
\def\uuug{{\n ^2 u+du\otimes du-\frac {|\n u|^2} 2 g+S_{g}}}
\def\sku{\sk\left(\uuu\right)}
\def\qed{\cqfd}
\def\vvv{{\frac{\n ^2 v} v -\frac {|\n v|^2} {2v^2} g_0+S_{g_0}}}
\def\vvs{{\frac{\n ^2 \tilde v} {\tilde v}
 -\frac {|\n \tilde v|^2} {2\tilde v^2} g_{\S^n}+S_{g_{\S^n}}}}
\def\skv{\sk\left(\vvv\right)}
\def\tr{\hbox{tr}}
\def\pO{\partial \Omega}
\def\dist{\hbox{dist}}
\def\RR{{\mathbb R}}
\def\R{{\mathbb R}}
\def\C{\Bbb C}
\def\B{\Bbb B}
\def\N{\Bbb N}
\def\Q{\Bbb Q}
\def\Z{\Bbb Z}
\def\PP{\Bbb P}
\def\EE{\Bbb E}
\def\G{\Bbb G}
\def\H{\mathbb H}
\def\S{{\mathbb S}}
\def\lcf{{locally conformally flat} }
\def\Lam{\Lambda}
\def\uuul{{\n ^2 u_i+du_i\otimes du_i-\frac {|\n u_i|^2} 2 g_0+S_{g_0}}}
\def\uuult{{\n ^2 \tilde u_i+d\tilde u_i\otimes d\tilde u_i-\frac
{|\n \tilde u_i|^2} 2 g_0+S_{g_0}}}
\def\uuulh{{\n ^2 \hat u_i+d\hat u_i\otimes d\hat u_i-
\frac {|\n \hat u_i|^2} 2 g_0+S_{g_0}}}
\def\uuug{{\n ^2 u+du\otimes du-\frac {|\n u|^2} 2 g+S_{g}}}
\date{}
\title[Geometric inequalities]
{Geometric inequalities on
locally conformally flat manifolds}
\author{Pengfei Guan}
\address{Department of Mathematics\\
 McMaster University\\
Hamilton, Ont. L8S 4K1, Canada.\\
Fax: (905)522-0935 }
\email{guan@math.mcmaster.ca}
\thanks{Research of the first author was supported in part by
NSERC Grant OGP-0046732.}
\author{Guofang Wang}
\address{Max-Planck-Institute for Mathematics in
the Sciences\\ Inselstr. 22-26, 04103 Leipzig, Germany}
\email{gwang@mis.mpg.de}
\begin{abstract}
Through the study of some elliptic and
parabolic fully nonlinear PDEs, we establish conformal
versions of the quermassintegral inequality, the Sobolev
inequality and the Moser-Trudinger inequality
for the geometric quantities associated to
the Schouten tensor on \lcf manifolds. \end{abstract}
\maketitle

\section{Introduction}
In this paper, we are interested in certain global geometric
quantities associated to the Schouten tensor and their
relationship to each other in conformal geometry. For an oriented compact
Riemannian manifold $(M,g)$ of dimension $n>2$, there is a
sequence of geometric functionals arising naturally in conformal
geometry, which generalize the Yamabe functional, introduced by
Viaclovsky in \cite{Jeff1} as curvature integrals of the Schouten
tensor. If we let $Ric_g$ and $R_g$ denote the Ricci tensor and
the scalar curvature of $g$ respectively, the Schouten tensor can
be written as
\[ S_g=\frac 1{n-2}\left(Ric_g-\frac {R_g}{2(n-1)}\cdot g\right).\]
Let $\s_k$ be the $k$th elementary symmetric function, and let
\[\sk(g):=\sk (g^{-1}\cdot S_g)\]
be the {\it $\sk$-scalar curvature} of $g$. We are interested in
the following  functionals defined by
\begin{equation}\label{functional}
{ \mathcal F}_k(g)=vol(g)^{-\frac{n-2k} n}\int_M \sk(g)\, dg, \quad
k=0,1,...,n,
\end{equation}
where $dg$ is the volume form of $g$. For convenience, set $\s_0(g)=1$.

We note that ${\mathcal F}_0(g)$ is simply the volume of $g$ and
${\mathcal F}_1(g)$ is the  Yamabe functional (up to a constant multiple).
Let us fix a
background metric $g_0$ and denote by $[g_0]$  the conformal class of
$g_0$.
The standard Sobolev inequality states that \[{\mathcal F}_0(g)
\le C_{[g_0]}({\mathcal F}_1^{\frac n{n-1}}(g)), \quad \mbox{ for } g\in [g_0],
\] for some constant
$C_{[g_0]}$
depending only on $[g_0]$. Equality holds if and only if the
metric $g$ is a Yamabe minimizer.
 For $g\in [g_0]$, we can write $g=e^{-2u}g_0$.  Let
$S_{g_0}$ and $S_g$ be the corresponding Schouten tensors of $g_0$
and $g$ respectively. There is a transformation formula relating
them:
\begin{equation}\label{trans}S_g=\nabla^2 u +du\otimes du-\frac{|\nabla
u|^2}{2}g_0+S_{g_0}.\end{equation} From this point of view,
${\mathcal F}_k(g)$ in (\ref{functional}) resembles the $k$th
quermassintegral in the theory of convex bodies. The classical
quermasstintegral inequalities suggest that there may exist
similar inequalities between the geometric quantities ${\mathcal
F_k(g)}$.

The focus of this paper is the investigation of the precise
relationship of the geometric quantities defined in
(\ref{functional}). As we will see, there exists a complete system
of sharp inequalities for ${\mathcal F}_k$ on \lcf manifolds, and
the cases of equality will be characterized by certain ``extremal"
metrics, just as in the case of the Yamabe problem. The main
tool used to obtain these geometric inequalities is the theory of parabolic fully
nonlinear equations, and this is one of the special features of this
paper.

When $(M,g_0)$ is a \lcf manifold and $k\neq n/2$, it was proved
by Viaclovsky in \cite{Jeff1} that the critical points of
${\mathcal F}_k$ in $[g_0]$ are the metrics $g$ satisfying
the equation
\begin{equation}\label{critical}\s_k(g)=constant.\end{equation}
For $k=n/2$, ${ \mathcal F}_{\frac n2}(g)$ is a constant in the
conformal class (see \cite{Jeff1}). In this case, there is another
functional which was discovered in \cite{BV}. A similar functional
was also found in \cite{CY2}. The functional in \cite{BV} is
defined by
\begin{eqnarray}\label{E}
{\cal E_{n/2}}(g)&=&-\int_0^1\int_M \s_{n/2}(g_t)udg_tdt,
\end{eqnarray}
where  $u$ is the conformal factor of $g=e^{-2u}g_0$ and
$g_t=e^{-2tu}g_0$. Unlike $\F_k$, $\En$ depends on the choice of
the background metric $g_0$. However, its derivative $\n \En$ {\it
does not} depend on the choice of $g_0$. The critical points of
$\En$ correspond to the metrics $g$ satisfying (\ref{critical})
for $k=n/2$ on \lcf manifolds. These important variational
properties are unknown for general (i.e. not locally conformally flat)
manifolds for $k> 2$. For
this reason, we will restrict ourselves to \lcf manifolds in our
study of ${\mathcal F}_k$ in this paper. Another advantage of
\lcf manifolds is the availability of the fundamental result of
Schoen-Yau \cite{SY} on developing maps, which is crucial for us
in establishing the gradient estimates for the conformal flow
(\ref{flow}) we consider below.

\medskip


We will establish three types of inequalities depending on the
range of $k$. More precisely, a {\it Sobolev type} inequality
(\ref{Sob-ineq}) is established for any $k<\frac n2$ and a {\it
conformal quermassintegral type inequality} (\ref{cquer})   for
any $k\ge n/2$. And, for the exceptional case $k=n/2$,  we
establish a {\it Moser-Trudinger type inequality} (\ref{MT-ineq})
for $\En$.

\medskip
Before giving precise results, let us first recall some notation and
definitions. Let
\[\Gamma_k^+=\{\Lambda=(\l_1,\l_2,\cdots, \l_n)\in \R^n\,|\,
\sigma_j(\Lambda)>0, \forall\, j\le k\}.\]
Let $\overline{\Gamma}_k^+$ be the closure of $\Gamma_k^+$.
A $C^{2}$ metric $g$ is said to be in $\Gamma_k^+$
if the eigenvalues of $g^{-1}\cdot S_g(x)\in \Gamma_k^+$
 for any
$x\in M$ (see \cite{GVW}).  Such a metric is also called {\it admissible}.
Similarly, a $C^{1,1}$ metric is said to be in $\overline{\Gamma}_k^+$ if
 the eigenvalues of $g^{-1}\cdot S_g(x)\in \overline\Gamma_k^+$
 for any $x\in M$.
For convenience, we set $\s_0(A)=1$ and $\s_0(g)=1$.
We denote
\[{\mathcal C}_k([g_0])=\{g \in C^{4,\a}[g_0] | g \in \Gamma_k^+\},\]
where $[g_0]$ is the conformal class of $g_0$. Here $0<\a<1$.

Unlike the Yamabe equation ($k=1$ in (\ref{critical})), equation
(\ref{critical}) is fully nonlinear when $k>1$. Moreover, from work in \cite{CNS}, it is necessary to
assume $g\in \Gamma_k^+$ (or $\Gamma_k^-=-\Gamma_k^+$) for $k>1$
in order to make use of elliptic
theory. This is a very restrictive condition. It is proved in
\cite{GVW} that any \lcf $(M,g)$ with $g\in \Gamma_{\frac n2}^+$
is a space form. Nevertheless, for $k<\frac n2$, there are plenty
of other \lcf manifolds. For example, $\S^{1}\times \S^{n-1}$ with
the standard product metric is in $\Gamma_k^+$ for any $k<\frac
n2$, and $\H^p\times \S^{n-p}$ with the standard product metric is
in $\Gamma_k^+$ for certain range of $2\le k$ when $p$ is
sufficiently small compared to $n$. Furthermore, if $(M_1,g_1)$ and
$(M_2,g_2)$ are two $n$-dimensional \lcf manifolds  with
$g_i \in \Gamma_k^+$ ($i=1,2$) for some $1\le k<\frac n2$, then
one can construct a \lcf metric $g \in \Gamma_k^+$ on the
connected sum of them (see also \cite{GLW2}). This type of
construction yields a wide class of \lcf manifolds with metrics in
$\Gamma_k^+$ for each $k<\frac n2$.

\medskip
We now state our main results. Since any metric $g\in
\Gamma_{n/2}^+$ is conformally equivalent to a metric of constant
sectional curvature (see \cite{GVW}), we will assume, whenever we
consider (\ref{E}) in this paper, that $g_0$ has constant
sectional curvature.

\begin{thm}\label{inqs} Suppose that $(M,g_0)$ is a compact, oriented and connected
\lcf manifold with
$g_0 \in \Gamma_k^+$ smooth and
 $g\in {\mathcal C}_k$. Let $0\le l< k\le n$.
\begin{itemize}
\item[({\bf A}).] {\bf Sobolev type inequality:} If
$0\le l < k < \frac n2$, then there is a positive constant $C_S=C_S([g_0], n,k,l)$
depending only on  $n$, $k$, $l$ and the
conformal class $[g_0]$ such that
\begin{equation}\label{Sob-ineq} \left(\Fk(g)\right)^{\frac 1{n-2k}}
\ge C_S \left(\Fl(g)\right)^{\frac 1{n-2l}}.\end{equation} If we
normalize $\int_M \sigma_l(g)dg=1$, then equality holds in
(\ref{Sob-ineq}) if and only if
\begin{equation}\label{quot-eq}
\frac {\s_k(g)}{\s_l(g)}=C^{n-2k}_S .\end{equation} There exists a
metric $g_E\in {\mathcal C}_k$ attaining equality in
(\ref{Sob-ineq}). Furthermore,
\begin{equation}\label{best1}
C_S \le C_{S, k,l}(\S^n)=
{\binom{n}{k}}^{\frac{1}{n-2k}}{\binom{n}{l}}^{\frac{-1}{n-2l}}
\left(\frac{\omega_n^2}{2^n}\right)^{\frac{k-l}{(n-2k)(n-2l)}},\end{equation}
where $\omega_n$ is the volume of the standard sphere $\S^n$.
\item[({\bf B}).] {\bf Conformal quermassintegral type
inequality:} If $n/2\le k \le n$, $1\le l<k$, then
\begin{equation}\label{cquer}
\left(\Fk(g)\right)^{\frac 1 k}\le
{\binom{n}{k}}^{\frac 1k}{\binom{n}{l}}^{-\frac 1l} \left(\Fl(g)\right)^{\frac 1l}.\end{equation}
Equality in (\ref{cquer}) holds if and only if $(M, g)$ is
a spherical space form.
\item[({\bf C}).] {\bf Moser-Trudinger type inequality:}
If $k=n/2$, then
\begin{equation}\label{MT-ineq}
(n-2l){\cal E}_{n/2}(g) \ge C_{MT} \left\{
\log\int_M \s_l(g)dg-\log \int_M \s_l(g_0)dg_0\right\},\end{equation}
where
\[C_{MT}=\int_M \s_{n/2}(g_0)dg_0= \frac{\omega_n}{2^{\frac n2}}
\binom{n}{\frac n2}.\] Equality holds in (\ref{MT-ineq}) if and
only if $(M,g)$ is a space form. The above inequality is also true
for $l>k=n/2$, provided $g \in {\cal {C}}_l$.
\end{itemize}
\end{thm}

When $l=0$ and $k=1$, inequality (\ref{Sob-ineq}) is the standard
Sobolev inequality (e.g., see \cite{Aubin}). For $l=0$ and $1\le k
<\frac n2$, inequality (\ref{Sob-ineq}) provides the optimal control of the
$L^{n}$-norm of $e^{-u}$ in terms of $\cal F_k$. For these reasons,
we call (\ref{Sob-ineq})  a Sobolev type inequality. We were
informed by Professor Alice Chang that Maria Del Mar Gonzalez
\cite{Maria} obtained inequality (\ref{best1}) independently for
$l=0$. We suspect that (\ref{best1}) should be true on general
compact manifolds. Inequality (\ref{Sob-ineq}) plays a key role in
the proof of a main result in \cite{GLW}.

Inequality (\ref{cquer})
is reminiscent of the classical quermassintegral
inequality (e.g., see \cite{GW2} for a discussion), which provides one
of the motivations for this paper. In the case $n=4, k=2$ and $l=1$, inequality
(\ref{cquer}) was proved earlier by Gursky in \cite{Gursky} for
general $4$-dimensional manifolds. Some cases of the inequality
were also verified in \cite{GW2} and \cite{GVW} for \lcf
manifolds.

Inequality (\ref{MT-ineq}) is similar to the Moser-Trudinger
inequality on compact Riemannian surfaces (see \cite{Tru1},
\cite{Moser}, \cite{Onofri} and \cite{Hong}). When $l=0$,
(\ref{MT-ineq}) was proved by Brendle-Viaclovsky and Chang-Yang in
\cite{BV} and \cite{CY2} using a result in \cite{GW2} on a fully
nonlinear conformal flow, and was referred to by them as the
Moser-Trudinger inequality. We also
refer to \cite{Beckner} for a different form of the Moser-Trudinger
inequality in higher dimensions.

\medskip
The Yamabe problem is the problem of finding a metric in the
conformal class $[g_0]$ which minimizes the functional ${ \mathcal
F}_1(g)$. The final solution of the Yamabe problem by Aubin
\cite{Aubin} and Schoen \cite{Schoen} is one of the triumphs of
geometric analysis. There have been several important
generalizations of Yamabe type problems in other settings such
as the CR Yamabe problem studied by Jerison-Lee \cite{JL}.
Viaclovsky \cite{Jeff1} considered another type of Yamabe
problem: for $1\le k\le n$, find a metric $g\in [g_0]$
satisfying equation (\ref{critical}), or equivalently the following
partial differential equation for the function $u$ in local
orthonormal frames with respect to $g_0$:
\begin{equation}\label{critical1}\s_k\left(u_{ij}+u_i u_j-\frac{|\nabla u|^2}{2}\delta_{ij}
+S_{ij}\right)=c e^{-2ku},\end{equation} where $S_{ij}$  are the
entries of $S_{g_0}$, $u_i$ and $u_{ij}$ are the first and second
order covariant derivatives of $u$ with respect to the local
orthonormal frames. When $k>1$, equation (\ref{critical1}) is
fully nonlinear. There has been much recent activity surrounding
equation (\ref{critical1}). Viaclovsky \cite{Jeff1} obtained some
existence results for $k=n$ under certain geometric conditions.
Equation (\ref{critical}) was solved by Chang-Gursky-Yang
\cite{CGY1, CGY2} in the case $n=4, k=2$. Gursky-Viaclovsky
studied the equation for the negative cone case in \cite{GV2}, and
more recently for the lower dimensional cases $n=3,4$ in
\cite{GV1}. For the \lcf manifolds, equation (\ref{critical}) was
completely solved by Guan-Wang \cite{GW2} and Li-Li \cite{Yanyan},
and for more general conformally invariant equations by
Guan-Lin-Wang \cite{ GLW} and Li-Li \cite{LL3}. We remark that for
$k\ge \frac n2$, $(M,g_0)$ is conformally equivalent to a space
form by \cite{GVW}, so any solution of (\ref{critical}) is a
maximizer of $\cal F_k$. But, for $k<\frac n2$, the previous solutions
to equation (\ref{critical}) obtained in \cite{GW2, Yanyan} are
not necessarily the Yamabe type solutions. That is, they are not
necessarily minimizers of the functional ${\mathcal F}_k$.
Similarly, neither are the solutions of (\ref{quot-eq}) on a locally
conformally flat manifold obtained in \cite{GLW} and \cite{LL3}
necessarily minimizers of a corresponding functional.

One immediate consequence of Theorem \ref{inqs}  is the existence of a
minimizer (or maximizer) solution of the general conformal quotient
equation (\ref{quot-eq}). 
For simplicity of notation, we will denote $\frac{\s_k(A)}{\s_l(A)}$ by
$\frac{\s_k}{\s_l}(A)$. By the transformation formula
(\ref{trans}), equation (\ref{quot-eq}) can be expressed in the fully nonlinear
form:
\begin{equation}\label{quot-eq11}
\frac{\s_k}{\s_l}\left(u_{ij}+u_i u_j-\frac{|\nabla
u|^2}{2}\delta_{ij} +S_{ij}\right)=c e^{-2(k-l)u}.\end{equation}

\begin{cor}\label{corex} For any $0\le l<k$, there is an extremal metric
$g_E\in {\cal C}_k$ satisfying equation (\ref{quot-eq11}) which is a global
minimizer (maximizer) of $\Fk$ in ${\cal C}_k$ with $\int_M
\sigma_l(g_E)dg_E=1$ ($\En(g_E)=1$ if $l=n/2$) when $k<n/2$
($k>n/2$); when $k=n/2$, there is a global minimizer of
$\En$.\end{cor}

\medskip

We will prove Theorem \ref{inqs} by studying an associated flow
equation. The approach of using flows to study geometric
inequalities was previously considered in different contexts by
various authors in \cite{Andrews, Chou, Wangxj, TW}. In the
conformal context, a simpler flow was introduced in \cite{GW2}
\begin{equation}\label{flow0}
\left\{\ba{rcl}
\ds\vs \frac {d}{dt}g
&=&-(\log\sk(g)-\log r_k(g))\cdot g, \\
g(0)&=&g_0,\ea\right.
\end{equation}
where $r_k(g)$ is given by
\[ r_k(g)=  \exp \left(\frac 1{vol(g)} \int_M \log \sk (g)
\, dg\right).\]
The global existence, regularity and convergence for the flow
(\ref{flow0}) was proved there.
For the purpose of this paper, we need to deal with
the following general fully nonlinear flow:
 \begin{equation}\label{flow}\left\{\begin{array}{rcl}
\ds \frac d{dt} g&=&\ds\vs-\left(\log \frac{\s_k(g)}{\s_l(g)}-\log
r_{k,l}\right)\cdot g,\\
g(0)&=&g_0,\ea\right.
\end{equation}
where
\[r_{k,l}=\exp\left(\frac{\int \s_l(g) \log(\s_k(g)\s_l(g)^{-1}) dg}
{\int \s_l(g) dg}\right)\] is defined so that the flow
(\ref{flow}) preserves $\int \s_l(g)dg$ when $l \not = n/2$
and ${\cal E_{n/2}}$ when $l=n/2$.
We have the following result for the flow (\ref{flow}).

\begin{thm}\label{mainthm} For any $C^{4,\a}$ initial metric
$g_0\in\Gamma_k^+$, the flow (\ref{flow}) has a global
solution $g(t)\in {\mathcal C}_k$ for any $t>0$.
Moreover, there is $h \in \cal{C}_k$
satisfying equation (\ref{quot-eq}) such that,
\[ \lim_{t \to \infty} \|g(t)-h\|_{C^{4,\a}(M)}=0. \]
\end{thm}

Flows (\ref{flow0}) and (\ref{flow}) are fully nonlinear conformal
equations. The first step is to obtain a Harnack inequality (i.e.,
the {boundedness} of $|\n g(t)|$, see (\ref{Harnack}) below). This
can be done using a fundamental result of Schoen-Yau on developing
maps and the method of moving planes as in \cite{Ye}. This is
where the assumption of \lcf is used crucially. For the flow
(\ref{flow0}), since the volume of the evolved metric is preserved,
$C^0$ {boundedness} is a direct consequence of the Harnack
inequality (see \cite{GW2} and \cite{Ye}). However the flow
(\ref{flow}) may not preserve the volume. This is the major
difference between the flows (\ref{flow0}) and (\ref{flow}). To deal
with this  problem we first use the Harnack inequality
(\ref{Harnack}) to bound $|\n ^2g(t)|$. From there we obtain a
$C^0$ bound using our previous results \cite{GW2} and the local
estimates in \cite{GW1}. The preliminaries for this are developed
in Section 2.

\medskip

We will also need local estimates for equation (\ref{quot-eq}) in
the case $0\le l<k<\frac n2$. These estimates are needed to prove
Theorem 1.A.  The existence of such local estimates is a very
special feature of this type of conformal fully nonlinear
equation and is crucial for the study of the elliptic equation
(\ref{quot-eq}).

\medskip

This paper is organized as follows. In Section 2, we list some
basic facts regarding $\frac{\s_k}{\s_l}$ and the flow (\ref{flow}) and prove
Theorem  \ref{inqs} for the case $l=0$. In Section 3, we establish
some local estimates for equation (\ref{quot-eq}),
 which will be used in a crucial way to
obtain a positive lower bound for the geometric functionals.
In Section 4, we study the flow (\ref{flow}) and prove
 Theorem \ref{mainthm}. Theorem \ref{inqs}
for the case $l>0$  will be proved in Section 5. In Section 6, we
discuss the best constant, which  appears in part ({\bf A}) of
Theorem \ref{inqs}.
\medskip

\noindent{\it Notation:} In the rest of the paper, $u_i, u_{ij},
\cdots$ denote the covariant derivatives of the function $u$ with
respect to some local orthonormal frames of the background metric
$g_0$, unless it is otherwise indicated.

\medskip

\noindent {\it Acknowledgement:} Part of this work was carried out
while the second author was visiting Mathematics Department,
McMaster University, in April 2002. He would like to thank the
department for its hospitality. We would also like to thank Eric
Sawyer and the referees for valuable suggestions regarding the
presentation of this paper.

\section{Some basic facts}
Let $\L=(\lambda_1, \dots, \lambda_n) \in \mathbf{R}^n$. The $k$th
elementary symmetric function is defined as
\[\sk(\lambda_1, \dots, \lambda_n) = \sum_{i_1 < \dots < i_k}
\lambda_{i_1} \cdots \lambda_{i_k}.\]
A real symmetric $n \times n$ matrix $A$ is said to lie  in  $\Gamma_k^{+}$
if its eigenvalues lie in $\Gamma_k^{+}$.
Let $A_{ij}$ be the $\{i,j\}$-entry of an
$n \times n$ matrix. Then for $0 \leq k \leq n$, the
$k$th {\it Newton transformation} associated with $A$ is defined to be
\[T_k({ A}) = \sigma_k({A}) I - \sigma_{k-1}({ A}) A + \cdots +
(-1)^k A^k.
\]
We have
\[T_k({ A}) ^i_j= \frac{1}{k!} \delta^{i_1 \dots i_k i}_{j_1 \dots j_k j}
A_{i_1 j_1} \cdots A_{i_k j_k},\]
where $ \delta^{i_1 \dots i_k i}_{j_1 \dots j_k j}$ is the
generalized Kronecker delta symbol.  Here  we use the summation convention.
By definition,
\[\sk({A})=\frac{1}{k!} \delta^{i_1 \dots i_k }_{j_1 \dots j_k }
A_{i_1 j_1} \cdots A_{i_k j_k},\quad
T_{k-1}({ A}) ^i_j= \frac{\partial \sk({ A})}{\partial A_{ij}}.\]
For $0< l < k \le n$, let
\[\tilde T_{k-1,l-1}({ A})=
\frac {T_{k-1}({ A})}{\sk({ A})}-\frac {T_{l-1}({A})}{\s_l({
A})}.\] It is important to note that if $A\in \Gamma_k^+$, then
$\tilde T_{k-1,l-1}({ A}) \text{ is positive definite.}$ This
follows the fact that the function $\frac{\s_k}{\s_l}(A)$ is
monotonic in $\Gamma_k^+$, that is if $B-A$ is semi-positive, then
$\frac{\s_k}{\s_l}(A)\le \frac{\s_k}{\s_l}(B)$. Another important
fact is that $(\frac{\s_k}{\s_l}(A))^{\frac{1}{k-l}}$ is concave
in $\Gamma_k^+$ (e.g., see \cite{Tru}).  Let
$\L_i=(\l_1,\l_2,\cdots,\l_{i-1},\l_{i+1},\cdots,\l_n).$

\begin{lem}[Newton-MacLaurin Inequality \cite{HPL}] \begin{equation}
\label{NMI} l(n-k+1)\s_l(\L)\s_{k-1}(\L) \ge
k(n-l+1) \s_k(\L)\s_{l-1}(\L). \end{equation}
\end{lem}
\begin{lem}[Garding's Inequality]\label{gardingi}
Let $\L=(\l_1,\cdots,\l_n), \L_0=(\mu_1,\cdots,\mu_n) \in
\Gamma_k^+$, \[ F(\L)=\left(\frac {\s_k}{\s_l}(\L)\right)^{\frac
1{k-l}}.\] Then,
\begin{equation}\label{EE}\sum_i\{\frac {\s_{k-1}(\L_i)}{\s_k(\L)} -\frac
{\s_{l-1}(\L_i)}{\s_l(\L)}\}\mu_i\ge (k-l) \frac{F(\L_0)}
{F(\L)}.\end{equation}
\end{lem}
\noindent{\it Proof.} The main argument of the proof follows from
\cite{CNS}.
 {From} the concavity of $F$ in $\Gamma_k^+$,
  for $\L, \L_0=(\mu_1,\cdots,\mu_n)
\in \Gamma_k^+$ we have
\begin{eqnarray*} F(\L_0) &\le& F(\L)+\sum_l(\mu_i-\l_i)
\frac{\partial F(\L)}{\partial \l_i}\\
&=&F(\L) +\frac{1}{k-l}F(\L)\sum_i\{\frac {\s_{k-1}(\L_i)}{\s_k(\L)}
-\frac {\s_{l-1}(\L_i)}{\s_l(\L)}\}(\mu_i-\l_i)\\
&=&\frac{1}{k-l}F(\L)\sum_i\{\frac {\s_{k-1}(\L_i)}{\s_k(\L)}
-\frac {\s_{l-1}(\L_i)}{\s_l(\L)}\}\mu_i .
\end{eqnarray*}
In the last equality, we have used the fact that $F$ is
homogeneous of degree one. Then (\ref{EE}) follows. \qed


\begin{lem}\label{cor_a} A conformal class of metric $[g]$ with $[g]\cap\Gamma_k^+\not =
\emptyset$ does not have a $C^{1,1}$ metric $g_1\in \overline \Gamma_k^+$
 with
$\s_k(g_1)=0$, where $\overline \Gamma_k^+$ is the closure of $\Gamma_k^+$.
\end{lem}
\pr
By the assumption, there is a smooth admissible metric $g_0$ with
$\sk(g)>0$. Assume by contradiction that there is a $C^{1,1}$
metric $g_1$ with $\sk(g_1)=0$. Write $g_1=e^{-2u}g_0$, so $u$
satisfies
\begin{equation}\label{viscosity}
\sk \left(\uuu\right)=0.\end{equation}
Let
\[W=\left(\uuu\right)\quad \text{and} \quad
a_{ij}(W)=\frac{\partial \sk(W)}{\partial w_{ij}}.\] (Here, the
notation $a_{ij}(W)=\frac{\partial \sk(W)}{\partial w_{ij}} $ is
usually used by analysts. More intrinsically,
$a_{ij}=T_{k-1}(W)^i_j$, where $T_{k-1}(W)$ is the $(k-1)$th
Newton transformation defined above and $W$ is viewed as the
matrix $g_0^{-1}\cdot W$.) We may assume $u\le 1$ and $u(x_0)=1$
for some $x_0\in M$, since $u+c$ also satisfies (\ref{viscosity})
for any constant $c$. Let $v=e^{-u}-e^{-1}$, $h_t=t e^{-u}
+(1-t)e^{-1}$, $u_t=-\log h_t$ and $W_t=\n^2 u_t +du_t \otimes d
u_t -\frac {|\n u_t|^2} {2} g_0 +S_{g_0}$. As in \cite{Jeff4}, one
can check that (a) $W_t \in \Gamma_k^+$ and (b) $ (a_{ij}(W_t))$
positive definite (nonnegative definite, resp.)
 for all $0\le t <1$
($0\le t \le 1$ resp.). We have the following
\begin{equation}\label{dvi}
0>\sk(W_1)-\sk(W_0)=\sum_{ij}  (\int_0^1 \frac{a_{ij}(W_t)}{h^2_t}dt)\n^2_{ij}v
+\sum_l b^l(t,x) \n_l v +dv,
\end{equation}
for some bounded functions $d$ and  $b^l$, $l=1,..., n$. Since
$v\ge 0$ and $v(x_0)=0$, this is a contradiction to the strong
maximum principle. \qed

\bigskip

For any $0<k\le n$, let
 \begin{eqnarray}\label{tff}{\tFF}_{k}(g)=\begin{cases}
 \ds\vs\frac{1}{n-2k}\int_M \s_k(g) dg, & k \not =n/2,\\
\ds  \En(g), & k=n/2.\end{cases} \end{eqnarray}

\begin{lem} Flow (\ref{flow}) preserves $\tFF_l$, it also decreases
the functional $\tFF_{k}$. In fact, we have the following
evolution identities. The evolution equations for $\log
\frac{\s_k}{\s_l}$ and $\tFF_k$ are
\begin{equation}\label{kl-eq}
\frac d {dt}\log \frac{\s_k}{\s_l}(g) =\frac 12 {\rm tr} \{\tilde
T_{k-1,l-1}(S_g)\n_g^2 \log \frac{\s_k}{\s_l}(g)
\}+(k-l)(\log\frac{\s_k}{\s_l}(g)-\log r_{k,l})\end{equation} and
\begin{equation}\label{4}
\frac d {dt}\tFF_k(g)=-\frac 12 \int_M
\left(\frac{\s_k}{\s_l}(g)-r_{k,l}\right) \left(\log
\frac{\s_k}{\s_l}(g)-\log r_{k,l}\right)\s_l(g)dg.
\end{equation}
\end{lem}

\pr It is easy to verify $\frac d {dt}\tFF_l(g)=0$ from the
equation and the definition of $r_{k,l}$. The first equality
follows the same argument in \cite{GW2}. We verify the second
equality. When $k\neq \frac n2$,

\begin{eqnarray*}
\begin{array}{rcl}
\ds \frac d {dt}\tFF_k(g) &=&\ds\vs
\frac 12 \int_M\s_k(g) g^{-1}\frac d {dt} g \,dg\\
&=&\ds-\frac 12 \int_M \left(\frac{\s_k}{\s_l}(g)-r_{k,l}\right)
\left(\log \frac{\s_k}{\s_l}(g)-\log r_{k,l}\right)\s_l(g)\,dg.
\end{array}\end{eqnarray*}
The first equality was proven in \cite{Jeff1}, where the
assumption of locally conformally flat is used. By \cite{BV}, the
above also holds for $k=\frac n2$. \qed

\bigskip

In the rest of this section we use results obtained in \cite{GW2}
 for flow (\ref{flow0})
to establish Theorem \ref{inqs} in the case $l=0$.

\begin{pro}\label{GW} Let $(M, g_0)$ be a locally conformally
flat manifold with $g_0\in \Gamma_k^+$.
We have
\begin{itemize}
\item [(a).] When $k>n/2$, there is a constant
$C_Q=C_Q(n,k)>0$ depending only on $n$ and $k$
such that for any metric $g\in {\cal C}_k$.
\[\int_M \s_k(g)vol(g) \le C_Q vol(g)^{\frac{n-2k}n}.\]
\item [(b).] When $k<n/2$, there is a constant
$C_{S}=C_{S}([g_0],n)>0$ such that for any metric $g\in {\cal
C}_k$.
\[\int_M \s_k(g)vol(g) \ge C_S vol(g)^{\frac{n-2k}n}.\]
\item [(c).] If $k=n/2$ and $g_0$ is a metric of constant sectional
curvature, 
then for any $g\in {\cal C}_k$
\[{\cal E}_{n/2}(g) \ge\frac 1n C_{MT} (\log vol(g)-\log vol(g_0)),\]
where $C_{MT}=\int_M \s_{n/2}(g_0) dg_0$.
\end{itemize}
Moreover, in  cases (a) and (c) the equality holds if and only if
$g$ is a metric of constant sectional curvature.
\end{pro}

\pr (c) was already proved in \cite{BV}.

When $k>n/2$, from \cite{GVW} we know that $(M, g_0)$ is
conformally equivalent to a spherical space form. In this case, it
was proved in \cite{Jeff1} that any solution of (\ref{quot-eq})
for $l=0$ is of constant sectional curvature.  By Theorem 1 in
\cite{GW2} ($k>\frac n2$),
for any $g\in{\cal C}_k$ there is a metric
$g_e\in{\cal C}_k$ of
constant sectional curvature with $vol(g)=vol(g_e)$ and
\begin{equation}\label{100}
\tFF_k(g)\ge \tFF_k(g_e).\end{equation}
When $k>n/2$, (\ref{100}) implies
that
\[ vol(g)^{-\frac {n-2k} n} \int_M \s_k(g)  dg \le
vol(g_e)^{-\frac {n-2k} n} \int_M \s_k(g_e) dg_e,\]
and the equality holds if and only if $(M,g)$ is a space form.
It is clear that
\[vol(g_e)^{-\frac {n-2k} n} \int_M \s_k(g_e) dg_e\]
is a constant depending only on $n, k$. This proves (a).


It remains to prove (b).
For this case, we only need to prove that
\[\inf_{{\cal C}_k\cap\{ vol(g)=1\}} \Fk(g) =: \b_0>0.\]
Assume by contradiction that $\beta_0=0$.
By Theorem 1 in \cite{GW2}, we can find a sequence of
solutions $g_i=e^{-2u_i}g_0\in {\cal C}_k$ of (3) with $vol(g_i)=1$
and $\sk(g_i)={\b_i}$ such that $\lim_{i\to \infty}{\b_i}=0$.
$\sk(g_l)={\b_i}$ means
\begin{equation}\label{appl_1}
\s_k(\uuul)={\b_i} e^{-2k{u_i}}.\end{equation}
Consider the scaled metric $\tilde g_i= e^{-2\tilde u_i}g_0$
with $\tilde {u}_i={u_i} -\frac 1{2k} \log {\b_i}$,
which satisfies clearly that
\begin{equation}\label{appl_2}
\s_k\left(\uuult\right)= e^{-2k\tilde {u}_i}\end{equation} and
\[vol(\tilde g_i)= {\b_i}^{\frac n{2k}} \to 0 \quad \text{ as } i\to \infty.\]
By Corollary 1.2 in \cite{GW1}, we conclude that
\[ \tilde {u}_i \to +\infty \text { uniformly as } i\to \infty.\]
Hence ${m_i}:=\inf_M \tilde {u}_i \to +\infty$ as $i\to \infty$.
Now at the minimum point $x_i$ of $\tilde u_i$, by (\ref{appl_2}),
\[\s_k(S_{g_0}) \le \s_k\left(\uuult\right)= e^{-2km_i} \rightarrow 0.\]
This is a contradiction to the fact $g_0 \in \Gamma_k^+$.\qed

\section{Local estimates for conformal quotient equations}
In order to get a positive lower bound for the constant in ({\bf
A}) of Theorem \ref{inqs}, as in the proof of Proposition \ref{GW}
in the previous section, we prove some estimates for solutions of
\begin{equation}\label{quot-eq1}
\frac{\s_k}{\s_l}(g)=f,
\end{equation}
defined locally in any open subset of $M$ with a nonnegative $C^2$
function $f$. When $l=0$ local estimates were established in
\cite{GW1}. The same estimates hold for equation (\ref{quot-eq1})
when $0\le l<k\le n$. Since the proof for the special case
$(n-k+1)(n-l+1)>2(n+1)$ is much simpler and it is suffice for the
purpose of this paper, we only treat this case here. The local
estimates for equation (\ref{quot-eq1}) in general case $0\le
l<k\le n$ will appear elsewhere. We emphasize that, in this
section, we {\it do not} assume $(M,g_0)$ is \lcf.
\begin{thm}\label{thm-local}
Let $(M,g_0)$ be a smooth compact, $n$-dimensional manifold and let
$g=e^{-2u}g_0$ be a $C^{4}$ solution of (\ref{quot-eq1}) in $B_r$ with
\begin{equation}\label{res}
(n-k+1)(n-l+1)>2(n+1).\end{equation}
 Then there exist two positive constants $c_1$
  (depending only on  $\|g_0\|_{C^3(B_r)}$, $n, k, l$ and
  $\|f\|_{C^1(B_r)})$) and $c_2$
(depending only on  $\|g_0\|_{C^4(B_r)}$, $n, k, l$ and $\|f\|_{C^2(B_r)}$)
such that $\forall x \in B_{\frac r2}$,
\begin{equation}\label{local1}
 |\nabla u(x)|^2\le c_1 (1+e^{-2 \inf_{B_r} u})
\end{equation}
and \begin{equation}\label{local2}
|\nabla^2 u(x)|\le c_2 (1+e^{-2 \inf_{B_r} u}).
\end{equation}
In particular, if  $k<\frac n2$, condition (\ref{res}) holds, and
hence
 inequalities (\ref{local1}) and (\ref{local2}) are valid in this case.
\end{thm}

Theorem \ref{thm-local} has the following consequence.

\begin{cor}\label{newenergy}
Let $k,l$ as in Theorem \ref{thm-local}. There exists a constant
$\e_0>0$ such that for any sequence of $C^{4,\a}$ solutions $u_i$ of
(\ref{quot-eq1}) in $B_r$ with
\[
\int_{B_1} e^{-nu}dvol(g_0)\le \e_0,\] either
\begin{itemize}
\item[(1)] There is a subsequence $u_{i_l}$ uniformly converging to
$+\infty$ in any compact subset in $B_r$, or \item[(2)] There is a
subsequence $u_{i_l}$ converges strongly in $C^{4,
\alpha}_{loc}(B_r)$. If $f$ is smooth and
strictly positive in $B_r$ and $u_i$ is smooth,
then $u_{i_l}$ converges strongly in
$C^{m}_{loc}(B_r)$, $\forall m$.
\end{itemize}\end{cor}

The proof of Corollary \ref{newenergy} follows from the same lines
in the proof of Corollary 1.2 in \cite{GW1}, the same argument can
be traced back to Schoen \cite{Schoen1}, we will not repeat it
here.

\medskip

\noindent{\it Proof of Theorem \ref{thm-local}.} The proof in
\cite{GW1} works for (\ref{local2}) as
$(\frac{\s_k}{\s_l})^{\frac{1}{k-l}}$ is elliptic and concave in
$\Gamma_k^+$. We only need to prove (\ref{local1}).

We follow the same lines of proof in \cite{GW1} to prove the Theorem,
together with  Lemma \ref{lem-B} proven at the end of this section.
For convenince of the reader,  we will use
the same notations  as in \cite{GW1}.

Without loss of generality, we may assume $r=1$. Let
\[F=\left(\frac{\s_k}{\s_l}(W)\right)^{\frac 1{k-l}} \quad \hbox{
and }\quad F^{ij}=\frac{\partial F}{\partial w_{ij}},\] where
$W=(\n ^2 u+du\otimes du-\frac {|\n u|^2}2 g_0+S_{g_0})$. More intrinsically,
 \[F^{ij}=\frac1{k-l}\left(\frac{\s_k}{\s_l}(W)\right)^{\frac 1{k-l}-1}\frac
 1 {\s^2_l (W)} \left\{
 \s_l(W)T_{k-1}(W)^i_j-\s_k(W)T_{l-1}(W)^i_j\right\},\]
 where $T_j$ is the $j$th Newton transformation.
 Let $\rho$ be a test function $\rho\in C^\infty_0(B_1)$ such that
\begin{equation}\begin{array}{l}\label{8}
\vs \ds \rho  \ge  0,  \hbox{ in } B_1 \quad \hbox{ and }\quad
\vs\ds \rho  = 1, \hbox{ in } B_{1/2},\\
\vs\ds |\n \rho (x)| \le 2 b_0 \rho^{1/2}(x) \quad \hbox{ and }\quad
|\n^2 \rho|  \le   b_0,    \hbox{ in } B_1.
\end{array}\end{equation}
Here $b_0>1$  is a constant depending only on the background metric $g_0$.
 Set $H=\rho |\n u|^2$. It suffices to bound $\max_{B_1}H$. Let
$x_0\in B_1$ be a maximum point of $H$.  We may assume that $W$ is
diagonal at the point $x_0$ by choosing a suitable normal
coordinates around $x_0$. Set $\l_i=w_{ii}$ and
$\L=(\l_1,\l_2,\cdots,\l_n)$. In what follows, all computations
are given at the point $x_0$. We denote $\L_i$ the vector with
$i$th component deleted from $\L$ and
\[F^*=\frac 1{k-l}\left(\frac{\s_k}{\s_l}(\L)\right)^{\frac 1{k-l}-1}
\frac 1{\s_l^2(\L)}.\] At the point $x_0$, $(F^{ij})$ is diagonal
and
\begin{equation}\label{B2}
 F^{ii} = F^*
\{\s_l(\L)\s_{k-1}(\L_i)-\s_k(\L)\s_{l-1}(\L_i)\}.
\end{equation}

Since $W$ is diagonal at $x_0$, we have at $x_0$
\begin{equation}\label{add0}
    w_{ii}  =  u_{ii}+u^2_i-\frac12 |\n u|^2+S_{ii}, \quad
    u_{ij}  = -u_iu_j-S_{ij}, \quad \forall i\neq j,
\end{equation}
where $S_{ij}$ are entries of $S_{g_0}$.

We may assume that
\[H(x_0)\ge b_0^2A_0^2,\]
for some constant $A_0$ will be fixed later. We may also assume
that
\begin{equation}\label{add1.1} |S_{g_0}| (x_0)\le A_0^{-1}|\n
u|^2(x_0).\end{equation} Otherwise, we are done.
Since $x_0$ is the maximum point of $H$,
we can get, for any $i$,
\begin{equation}\label{add2.2}
\left|\sum_{l=1}^nu_{il}u_l\right|(x_0)\le \frac{|\n
u|^3}{A_0}(x_0).\end{equation}

By applying maximum principle and
following the same deduction in \cite{GW1} (formulas (2.18) and
(2.20) there), we have
\begin{equation}\label{add100} \begin{array}{rcl}
0 &\ge &\ds\vs
 F^{ij}H_{ij}=F^{ij}
\left\{\left(-2\frac{\rho_i\rho_j}\rho+\rho_{ij}\right)|\n u|^2
+2\rho u_{lij}u_l+2\rho u_{il}u_{jl}\right\}\\
&\ge & \ds
2\sum_{i,j,l}F^{ij}\rho u_{il}u_{jl}
+\sum_{j\ge 1} F^{jj}\left\{-10n  b_0^2|\n u|^2-Ce^{-\inf_{B_1}u}|\n u|^2
-\frac{(n+2)^2}{A_0}\rho |\n u|^4\right\}
\end{array}\end{equation} where $C$ is a constant
depending only on $n, k, l$, $\|f\|_{C^1(B_1)}$ and
$\|g\|_{C^3(B_1)}$. Set $\tilde u_{ij}=u_{ij}+S_{ij}$. By
(\ref{add1.1}), we have
\begin{eqnarray*}
\ds\vs \sum_{i,j,l}F^{ij}u_{il}u_{jl} \ge \frac12
\sum_{i,l}F^{ii}\tilde u^2_{il} -\frac n {A^2_0} |\n
u|^4\sum_{i}F^{ii}.\end{eqnarray*} This, together with (\ref{add100}),
implies
\begin{eqnarray*}
\sum_{i,l}F^{ii}\tilde u^2_{il}\le
\sum_{j\ge 1} F^{jj}\left\{10n  b_0^2|\n u|^2+Ce^{-\inf_{B_1}u}|\n u|^2
+\frac{2(n+2)^2}{A_0}\rho |\n u|^4\right\}.
\end{eqnarray*}
Hence, the needed bound
$H(x_0)\le c_1(1+e^{-2\inf_{B_1}u})$ can be reduced to the
verification of the following:

{\bf Claim:} There is a constant $A_0$ (depending only on $k$ and
$n$) 
such that
\begin{equation}\label{B1}
\sum_{i,j}F^{ii}\tilde u^2_{ij} \ge A_0^{-\frac 58}\sum_i
F^{ii}|\n u|^4.
\end{equation}

This claim is just Claim 2.5 in \cite{GW1}. As in \cite{GW1}, the
verification of the Claim
is a crucial step to prove the Theorem.

Set $\d_0= A_0^{-1/4}<0.1$. We divide the set $I=\{1,2,\cdots,
n\}$ as in \cite{GW1} into two parts:
\[I_1=\{i\in I\,|\,u_i^2 \ge \d_0 |\n u|^2\} \quad \hbox{ and } \quad
I_2=\{i\in I\,|\,u_i^2 < \d_0 |\n u|^2\}.\] Clearly, $I_1$ is
non-empty.

\medskip

{\noindent \bf Case 1.}  There is $j_0$ satisfying
\begin{equation}\label{B4}
\tilde u_{jj}^2 \le \d_0^2 |\n u|^4 \quad \hbox{ and } \quad
u^2_{j} <\d_0|\n u|^2.\end{equation} We may assume that $j_0=n$.
We have $|w_{nn}+\frac {|\n u|^2} 2|= |\tilde u_{nn}+ u_n^2|<2\d_0
|\n u|^2$ by (\ref{B4}). From (\ref{add2.2}) and (\ref{add0}), we
have for any $i\in I$
\[\left|u_i(u_{ii}-(|\n u|^2-u_i^2))-\sum_lS_{il}u_l\right|=
\left|\sum_lu_{il}u_l\right|\le \d_0^4|\n u|^3.\] This, together
with (\ref{add1.1}), implies
\begin{equation}\label{add3}
  |u_i(u_{ii}-(|\n u|^2-u_i^2))|\le 2\d_0^4|\n u|^3,
\end{equation}
which, in turn, implies \[\begin{array}{rcl} \left|w_{ii}-\frac
{|\n u|^2} 2\right|&=& \left|u_{ii}+u_i^2-|\n u|^2+S_{ii} \right|
\le 3\d_0^2 |\n u|^2,\end{array}\] for any $i\in I_1$. 
Using these inequalities, we can repeat the derivation of equation
(2.38) in \cite{GW1} to obtain
\begin{equation}\label{B3}
\ds\vs \sum_{i,l} F^{ii} \tilde u^2_{il}\ge
\tilde F^1\frac {|\n u|^4} 4-\tilde F^1 |\n u|^4+
F^{nn}\frac {|\n u|^4} 4+(1-32 \d^2_0) \frac{|\n u|^4} 4\sum_{i} F^{ii},
\end{equation}
where
$\tilde F^1=\max_{i\in I_1} F^{ii}$.

 Recall that $I_1$ is necessarily non-empty. We may assume that
$i_0\in I_1$ with $\tilde F^1=F^{i_0i_0}$. Without loss of
generality, we assume that $i_0=1$. {From} above, we know $w_{11}
\ge w_{nn}$. By Lemma \ref{lem-B} below, we have $F^{nn} \ge
\tilde F^{11}$ and
 $\tilde F^{11}\le \frac 1{2+c(n,k,l)}\sum_l F^{ll}$, because $w_{11}>0$.
 Hence,
 (\ref{B3}) implies that
 \[ \sum_{i,l} F^{ii} \tilde u^2_{il}
\ge
 \{\frac 14-\frac 1 {2(2+c(n,k,l))}-32 \d_0^2\} \sum_i F^{ii}|\n u|^4.\]
The inequality (\ref{B1}) then follows for this case by 
choosing $A_0\ge \left(\frac{128(2+c(n,k,l))}{c(n,k,l)}\right)^2.$

\medskip

{\noindent \bf Case 2.}
There is no $j\in I$ satisfying (\ref{B4}).

\medskip

We may assume that there is $i_0$ such that $\tilde
u_{i_0i_0}^2\le \d_0^2|\n u|^4$, otherwise the claim is
automatically true. Assume $i_0=1$. As in Case 4 in \cite{GW1}, we
have $u_1^2\ge (1-2\d_0)|\n u|^2$ and $\tilde u^2_{jj} + u_j^2(|\n
u|^2-u_j^2)\ge \d^2_0 |\n u|^4$ for $j>1$. From equation (2.50) in
\cite{GW1} and Lemma \ref{lem-B}, we have
\begin{eqnarray*}
\sum_{i,l}F^{ii}\tilde u_{il}^2 & =& \sum_{i}F^{ii}\left\{\tilde u_{ii}^2 +u_i^2(|\n u|^2-u_i^2)
\right\}\\
& \ge &  \sum_{i\ge 2}F^{ii}(\tilde u^2_{ii}+u_i^2(|\n u|^2-u_i^2))\\
&\ge& \d^2_{0} |\n u|^4\sum_{i\ge 2} F^{ii} \ge \frac 12
\d^2_{0} |\n u|^4\sum_{i\ge 1} F^{ii}
 \end{eqnarray*}
The {\bf Claim} is verified, so the local gradient estimate ({\ref{local1})
is proved.\qed


\begin{lem}\label{lem-B}
 If $\l_i \le \l_j$, then $F^{ii} \ge F^{jj}$. If $(n-k+1)(n-l+1)>2(n+1)$
 and $\l_1 >0$, then there is a positive constant $c(n,k,l)$ such that
 \[\sum_i F^{ii} \ge (2+c(n,k,l))F^{11}.\]
\end{lem}
\pr  The first statement follows from the monotonicity of
$\s_{l-1}$ and $\frac{\s_{k-1}}{\s_{l-1}}$, since
\[ F^{ii} =F^*\s_{l-1}(\L_i)\{\s_l(\L)\frac {\s_{k-1}}{\s_{l-1}}(\L_i)
-\s_k(\L)\}.\] Similarly, we have
\[\begin{array}{rcl}
\ds F^{11}& =& \ds\vs F^*\s_{l-1}(\L_1)\{\s_l(\L)\frac
{\s_{k-1}}{\s_{l-1}}(\L_1)
-\s_k(\L)\}\\
&\le &\ds F^*\s_{l-1}(\L)\{\s_l(\L)\frac {\s_{k-1}}{\s_{l-1}}(\L)
-\s_k(\L)\}. \end{array}\]
{{From}} the Newton-MacLaurin inequality (\ref{NMI}),
\[q_0:= (n-k+1-\a_0) \s_l(\L)\s_{k-1}(\L)-
(n-l+1-\a_0) \s_k(\L)\s_{l-1}(\L)\ge 0,\] where $\a_0=\frac
{(n-k+1)(n-l+1)}{n+1}$. Hence, we have
\begin{eqnarray*}
 \ds\sum_iF^{ii}& =&
 \vs F^*\{(n-k+1)\s_l(\L)\s_{k-1}(\L)-(n-l+1) \s_k(\L)\s_{l-1}(\L)\\
 &=& \ds\vs F^*\{q_0+ \a_0 [s_l(\L)\s_{k-1}(\L)-\s_k(\L)\s_{l-1}(\L)]\} \\
 & \ge & \a_0 F^{11}
 .\end{eqnarray*}
The Lemma follows with $c(n,k,l)=2-\a_0>0$.\qed

\section{The global convergence of flow (\ref{flow})}

We treat flow (\ref{flow}) in this section. If $g=e^{-2u}\cdot g_0$,
one may compute that
\[\sk(g)=e^{2ku}\sk\left(\uuu\right).\]
Equation (\ref{flow}) can be written in the
following form
\begin{equation}\label{uflow}
\left\{\ba{rcl}
\ds \vs 2\frac{\ds du}{\ds dt} & = &  \ds \log \frac{\sk}{\sigma_l}
\left(\uuu\right)+2(k-l)u -\log r_{k,l}(g)\\
u(0) & = & u_0.
\ea\right.\end{equation}

In these equations and below, the norms and covariant derivatives
are with respect to the background metric $g_0$.

\medskip

Because $g_0 \in \Gamma_k^+$, the highest order term on the right
hand side is uniformly elliptic. Consequently, the short time
existence of flow (\ref{flow}) follows from the standard parabolic
theory. We want to prove the long time existence and convergence.
Let
\[T^*=\sup\{T_0>0\,|\, \hbox{(\ref{flow}) exists in } [0,T_0]
\hbox{ and } g(t)\in {\mathcal C}_k \hbox{ for } t\in [0.T_0]\}.\]
 For any
$T<T^*$, we will establish a $C^2$ bound for the conformal factor
$u$, which is independent of $T$.
As in \cite{GW2}, we have a Harnack inequality
\begin{equation}\label{Harnack}
  |\n u|\le c,
\end{equation}
for some positive constant $c$ independent of $T$,
 see a complete proof in
\cite{Ye}. But, as mentioned above, we can not deduce $C^0$
{boundedness} as in \cite{GW2} because the  flow (\ref{flow}) may
not preserve the volume of the evolved metric $g(t)$. In order to
to get $C^1$ bound we first  use (\ref{Harnack}) to bound $|\n
u|^2$.

\begin{pro}\label{pro1} There is a constant $C>0$ independent of
$T$ such that
\[|\n^2 u|\le C.\]
\end{pro}
\proof Set \[F=\log \frac{\sigma_k}{\sigma_l}\left(\uuu\right).\]
By equation (\ref{uflow}), $F=2u_t-2(k-l)u-\log r_{k,l}.$ We only
need to consider the case $k>1$. From the property of $\Gamma_2$,
\[|\l_i|\le \max\{1,(n-2)\}\sum_{j=1}^n \l_j.\]
Therefore we only need to give a upper bound of $\D u$ which
dominates all other second order derivatives.

Consider $G= \D u
+m|\n u|^2$ on $M\times [0,T]$, where $m$ is a large constant
which will be fixed later. Note that $\D u$ is
analyst's Laplacian. Without loss of generality, we may assume
that the maximum of $G$ on $M\times [0, T]$ achieves at a point
$(x_0, t_0)\in M\times (0,T]$ and $G(x_0,t_0) \ge 1$. We may
assume that at $(x_0,t_0)$
\begin{equation}\label{c0}2\s_1(W)\ge G \ge \frac 12 \s_1(W),\end{equation}
 where $W=\uuu.$
Consider everything in a small neighborhood near $x_0$. We may consider
$W$ as a matrix with entry
$w_{ij}=u_{ij}+u_iu_j-\frac 12|\n u|^2 \d_{ij}+S(g_0)_{ij}.$
In the rest of the proof, $c$ denotes a positive constant
independent of $T$, which may vary from line to line.

Since $G$ achieves its maximum at $(x_0,t_0)$, we have at this
point
\begin{equation}\label{c2}
  G_t=\sum_l (u_{llt}+2mu_{lt}u_l) \ge 0,
\end{equation}
and
\begin{equation}\label{c3}
  G_i=\sum_l(u_{lli}+2mu_{li}u_l)=0, \quad \forall\,\, i.
\end{equation}
(\ref{c3}) and (\ref{Harnack}) imply that at $(x_0, t_0)$
\begin{equation}\label{c4}
|\sum_l u_{lli}|\le cG.
\end{equation}
By the Harnack inequality (\ref{Harnack}) and the fact
$|u_{ij}|\le G$, we may assume that
\begin{equation}\label{c5}
  |u_{lij}-u_{ijl}|<c \quad \text{ and } \quad
|u_{ijkl}-u_{ijlk}|< cG,
\end{equation}
where $c$ is a constant depending only on the background metric
$g_0$. We may assume by choosing appropriate local orthonormal
frames that the matrix $(w_{ij})$ at $(x_0,t_0)$ is diagonal. At
the maximum point, $G_{ij}$ is non-positive definite. Set
$F^{ij}=\frac {\partial F}{\partial w_{ij}}.$ Since
$g(t)=e^{-2u(t)}g_0\in \Gamma_k^+$, we know that the matrix
$(F^{ij})$ is positive. Hence in view of (\ref{c2})-(\ref{c5}) and
the concavity of $F$ we have
\begin{equation}\label{c7}\begin{array}{rcl}
\quad 0&\ge& \vs\ds \sum_{i,j} F^{ij}G_{ij} =\sum
F^{ij}(u_{llij}+2mu_{li}u_{lj}+2mu_{lij}u_l)\\
 &\ge &  \vs\ds \sum_{i,j,l} F^{ij}(u_{ijll} +2mu_{li}u_{lj}
 +2mu_{ijl}u_l)-c\sum_i F^{ii}G\\
 &=&
 \vs \ds  -c\sum_i F^{ii}G +\sum_{i,j,l}  F^{ij}\{w_{ijll}-
 (u_iu_j-\frac 12 |\n u|^2\d_{ij}+S(g_0)_{ij})_{ll}\\
 &&\vs\ds + 2mu_{li}u_{lj}+ 2mw_{ijl}u_l
  -2mu_l(u_iu_j-\frac 12 |\n u|^2\d_{ij}+S(g_0)_{ij})_l\}\\
& \ge &  \vs\ds
 \D F +2m \sum_l F_l u_l+\sum_i F^{ii} u^2_{jl}
 +2(m-1)\sum_{i,l} F^{ii} u^2_{li}-c\sum_i F^{ii}G\\
 & \ge &\ds \D F +2m \sum_lF_l u_l+\frac 1n G^2\sum_i F^{ii}+2(m-1)
 \sum_{i,l} F^{ii} u^2_{li}-c\sum_i F^{ii}G.
 \end{array}\end{equation}
{From} equation (\ref{uflow}), $ F=2u_t-2(k-l)u-\log r_{k,l}(g).$
In view of (\ref{c2}) and (\ref{c3}), (\ref{c7}) yields
\begin{equation}\label{c8}
\begin{array}{rcl}
0 & \ge & \vs\ds -2(k-l) G+\frac 1n\sum_{i}F^{ii}G^2+2(m-1) \sum_i
 F^{ii}u_{ii}
-c \sum_i F^{ii} G\\
&\ge & \vs \ds -2(k-l) \D u+\sum_{i} F^{ii} G^2+2(m-1) \sum
F^{ii}u^2_{ii}
-c \sum F^{ii} G\\
&\ge & \ds \{-2(k-l) G+2(m-1) \sum F^{ii}u^2_{ii}\}+\frac 1n
\sum_{i}F^{ii}(G^2-c G).
\end{array}\end{equation}
The Proposition follows
from (\ref{Harnack}), (\ref{c8}) and
Lemma \ref{claim40} below. \qed

\begin{lem}\label{claim40} There is a large constant $m>0$ such that
\begin{equation}\label{c9}
  \frac 1{2n}G^2\sum_i F^{ii}+2(m-1)\sum_iF^{ii}w_{ii}^2\ge
  2(k-l)G.
\end{equation}
\end{lem}

\pr
It is easy to check, from the Newton-MacLaurin inequality
(\ref{NMI}), that
\begin{equation}\label{c10}
\begin{array}{rcl}
 \vs\ds\sum F^{ii}w_{ii}^2&=&\ds\frac{\s_1(W)\s_k(W)-(k+1)\s_{k+1}(W)}
 {\s_k(W)}-\frac{\s_1(W)\s_l(W)-(l+1)\s_{l+1}(W)}
 {\s_l(W)} \\
 &=& \vs\ds
 (l+1)\frac{\s_{l+1}}{\s_l}(W)-(k+1)\frac{\s_{k+1}}{\s_k}(W)
 \ge \ds c_{n,k,l}\frac {\s_{l+1}}{\s_l}(W),
 \end{array}
\end{equation}
and
\begin{equation}\label{c11}
 \sum_iF^{ii} =\vs\ds(n-k+1)\frac {\s_{k-1}}{\s_k}(W)-
  (n-l+1)\frac{\s_{l-1}}{\s_l}(W)
\ge \frac{\tilde c_{n,k,l}}{\sigma_1(W)},
\end{equation}
where $c_{n,k,l}$ and $\tilde c_{n,k,l}$ are two positive constant
depending only on $n, k$ and $l$. {From} these two facts, we can
prove the claim as follows. First, if
\[\frac {\tilde c_{n,k,l}}{4n}\frac
{\s_1(W)\s_{k-1}(W)}{\s_k(W)}\ge 4(k-l),\] then the claim follows
from (\ref{c11}) and (\ref{c0}). Hence we may assume that
\[\frac{\s_1(W)\s_{k-1}(W)}{\s_k(W)}
\le c^*_{n,k,l},\] for some positive
constant $c^*_{n,k,l}$ depending only on $n, k$ and $l$. Together
with  the Newton-MacLaurin inequality, it implies
\[\frac {\s_{l+1}(W)}{\s_{l}(W)}\ge \hat c_{n,k,l}\frac
{\s_{k}(W)}{\s_{k-1}(W)} \ge \hat c_{n,k,l} c^*_{n,k,l}\s_1(w),\]
which, in turn, together with (\ref{c10}) implies
\[\sum_iF^{ii}w_{ii}^2\ge c_{n,k,l}\frac {\s_{l+1}(W)}{\s_l(W)}\ge
c^1_{n,k,l}\s_1(W)\ge \frac 12c^1_{n,k,l} G.\] Hence, if we choose
$m$ large, then the lemma is true.  \qed

Now we can prove the $C^0$ {boundedness} (and hence $C^2$ {boundedness}).
\begin{pro}\label{pro2}
Let $g=e^{-2u}g_0$ be a solution of flow (\ref{flow})
with $\s_k(g(t)) \in \Gamma^+_k$ on $M\times [0,T^*)$. Then there
is a constant $c>0$ depending only on $v_0$, $g_0$, $k$ and $n$ (independent of $T^*$) such that
\begin{equation}\label{b20}
\|u(t)\|_{C^2}\le c,\quad \forall\, t \in [0, T^*).
\end{equation}
\end{pro}
\pr  We only need to show the {boundedness} of $|u|$.
First we consider the case $l\neq n/2$.
By  Proposition \ref{pro1} and the preservation of $\int \s_l(g)
dg$, we have
\begin{equation}
\label{c21}\begin{array}{rcl} c_l& =& \ds\vs\int_M
e^{(2l-n)u}\s_l\left(\uuu\right)dg_0\\
&\le &
c_1\int_Me^{(2l-n)u}dg_0.\end{array}
\end{equation} If $l<n/2$, then
(\ref{c21}), together with (\ref{Harnack}), implies that $u<c$ for
some constant $c>0$. On the other hand, in this case Proposition \ref{GW} gives
\[vol(g) \le C \left(\int_M \s_l(g) dg\right)^{\frac n{n-2l}}=c_0C,\]
which, together with (\ref{Harnack}) implies $u>c_1$, hence
$|u|\le C$ in this case.

If $l>n/2$, (\ref{c21}) gives a lower bound of $u$.
Suppose there is no upper bound, we have a sequence of
$u$, with $\nabla u$ and $\nabla^2 u$ bounded, but
$\sup u$ goes to infinity (so does $\inf u$). Set
$v=u- \inf u$, so $v$ is bounded and so is the
$C^2$ norm of $v$. But, for $\tilde g =e^{-2v}g_0$,
we get $\tilde F_l(\tilde g)$ tends to $0$. Take
a subsequence, we get $\s_l(e^{-2v^*}g_0)=0$ with $v^*$
in $C^{1,1}\cap \overline {\Gamma}_k^+$. This is a contradiction
to Lemma \ref{cor_a}.

Then we consider the case $l=n/2$. In this case, $\En(g)$ is constant.
First it is easy to check that $g_t=e^{-2tu}g_0 \in \Gamma_{\frac n2}^+$
when $0\le t \le 1$ (using the fact
$(1,\cdots,1,-1) \in \overline \Gamma_{\frac n2}^+$ when $n$ even).
In particular,
$\sigma_{\frac n2}(g_t)>0$ for $t>0$. {From} the expression of $\En(g)$,
\[-\sup(u)\int_M \sigma_{\frac n2}(g)dg \le
\En(g) \le - \inf(u) \int_M \sigma_{\frac n2}(g)dg.\]
Since
\[\int_M \sigma_{\frac n2}(g)dg=\int_M \sigma_{\frac n2}(g_0)dg_0,\]
So we have
\[-\sup(u)\int_M \sigma_{\frac n2}(g_0)dg_0 \le
\En(g) \le - \inf(u) \int_M \sigma_{\frac n2}(g_0)dg_0.\] Thus,
$\inf(u)$ is bounded from above and $\sup(u)$ is bounded from
below. By (\ref{Harnack}) again, $u$ is bounded from above and
away from $0$. Now we have proved the  {boundedness} of $|u|$ in
all cases. Combining that with (\ref{Harnack}) and Proposition
\ref{pro1} gives  a $C^2$ bound $u$ independent of $T$. \qed

\begin{pro}\label{pro3}
There is a constant $C_0>0$ independent of $T$ such that
\[\frac{\s_k}{\s_l}(g)(t) \ge C_0, \quad \text{ for } t \in [0, \infty).\]
\end{pro}
\pr We follow the similar argument in \cite{GW2}. Here we will
make use of Lemma \ref{gardingi}. We consider $H=\log
\frac{\s_k}{\s_l}(g)-e^{-u}$ on $M\times [0,T]$ for any $T<T^*$.
 {From} (\ref{flow}) and (\ref{kl-eq}) we have
\[\begin{array}{rcl}
\ds \frac{dH}{dt} &=& \ds\vs \frac 12\tr \{\tilde T_{k-1,
l-1}(S_g)\n_g^2 \log \frac{\s_k}{\s_l}(g)\}+
(k-l+\frac 12e^{-u})\left(\log\frac{\s_k}{\s_l}(g)-\log r_{k, l}(g)\right)\\
&=& \ds\vs  \frac 12\tr \{\tilde T_{k-1, l-1}(S_g)\n_g^2
(H+e^{-u})\}+ (k-l+\frac
12e^{-u})\left(\log\frac{\s_k}{\s_l}(g)-\log r_{k, l}(g)\right).
\end{array}
\]
Without loss of generality, we may assume that the minimum of $H$
in $M\times [0, T]$  achieves at $(x_0, t_0)\in M\times (0, T]$.
Let $H_{j}$ and $H_{ij}$ are the first and second derivatives with
respect to the background metric $g_0$. At this point,  we  have
$\frac {dH}{dt} \le 0$, $0= H_l = \sum_{ij}F^{ij}
w_{ijl}+e^{-u}u_l$ for all $l$, and $(H_{ij})$  is non-negative
definite. Also we have $(F^{ij})$ is positive definite and
\[\sum_{i,j} F^{ij}w_{ij}=\frac 1 {\s_k(g)}
\frac {\partial \s_k(g)}{\partial w_{ij} } w_{ij}
-\frac 1 {\s_l(g)}
\frac {\partial \s_l(g)}{\partial w_{ij} } w_{ij}=k-l.\]
Recall that in local coordinates
$\tilde T^{ij}_{k-1, l-1}(S_g)=F^{ij}$
and
\[\sum_{i,j} F^{ij} (\n^2_g)_{ij} H
= \ds \sum_{i,j}F^{ij}(H_{ij}+ u_iH_j+u_jH_j-\sum_l u_lH_l \d_{ij}).\]
It follows that at the point,
\begin{equation}
\label{add1}\begin{array}{rcl}
0& \ge &  \vs\ds H_t-\frac 12 \sum_{i,j}F^{ij} H_{ij} \\
&=&\vs\ds \frac 12\tr \{\tilde T_{k-1, l-1}(S_g)\n_g^2 e^{-u}\}+
(k-l+\frac 12 e^{-u})\left(\log\frac{\s_k}{\s_l}(g)-\log r_{k,l}(g)\right)\\
&=&\vs\ds \frac 12 \sum_{i,j} F^{ij}
\{(e^{-u})_{ij}+u_i(e^{-u})_{j}+u_j(e^{-u})_{i}-u_l(e^{-u})_l\d_{ij}\}\\
&& \ds\vs +(k-l+\frac 12 e^{-u})\left(\log\frac{\s_k}{\s_l}(g)-\log r_{k,l}(g)\right)\\
&=&\vs\ds \frac{e^{-u}}{2}\sum_{i,j} F^{ij}\{-u_{ij}-u_iu_j+|\n u|^2\d_{ij}\}
+(k-l+\frac 12e^{-u})\left(\log\frac{\s_k}{\s_l}(g)-\log r_{k,l}(g)\right)\\
&=&\vs\ds\frac{e^{-u}}{2}\sum_{i,j}  F^{ij}\{-w_{ij}+S_{ij}+\frac 12|\n u|^2\d_{ij}\}+
(k-l+\frac {e^{-u}}{2})\left(\log\frac{\s_k}{\s_l}(g)-\log r_{k,l}(g)\right)\\
&\ge & \vs\ds \frac{e^{-u}}{2} \sum_{i,j} F^{ij}\{-w_{ij}+S_{ij}\}+
(k-l+\frac 12e^{-u})\left(\log\frac{\s_k}{\s_l}(g)-\log r_{k,l}(g)\right)\\
&=&\ds \frac{e^{-u}}{2} \sum_{i,j} F^{ij}S_{ij} + (k-l+\frac
12e^{-u})\left(\log\frac{\s_k}{\s_l}(g)-\log r_{k,l}(g)\right)
-\frac{k-l}{2}e^{-u},
\end{array}\end{equation}
where $S_{ij}$ are the entries of $S(g_0)$.
Since $S(g_0) \in \Gamma_k^+$, by Lemma \ref{gardingi},
\begin{equation}
\label{add2} F^{ij}S_{ij}= \{\frac 1{\s_{k}(g)}\frac{\partial
\s_{k}(g)}{\partial w_{ij}} -\frac 1{\s_{l}(g)}\frac{\partial
\s_{l}(g)}{\partial w_{ij}}\} S_{ij} \ge
(k-l)e^{2u}\left(\frac{\s_k}{\s_l}(g_0)\right)^{\frac1{k-l}}
\left(\frac{\s_k}{\s_l}(g)\right)^{-\frac1{k-l}}.\end{equation} By
$C^2$ estimates, $\log r_{k,l}(g)$ is bounded from above, we have
\[ \begin{array}{rcl}
0&\ge & \ds\vs
\frac{(k-l)e^u}{2}\left(\frac{\s_k}{\s_l}(g_0)\right)^{\frac1{k-l}}
\left(\frac{\s_k}{\s_l}(g)\right)^{-\frac1{k-l}}\\
&&\ds\vs+(k-l+\frac 12e^{-u})\left(\log\frac{\s_k}{\s_l}(g)-\log
r_{k,l}(g)\right)-
\frac{k-l}{2}e^{-u}\\
&\ge & \ds c_1\left(\frac{\s_k}{\s_l}(g)\right)^{-\frac1{k-l}}+
c_2\log\frac{\s_k}{\s_l}(g)-c_3\end{array}\] for positive
constants $c_1$, $c_2$  and $c_3$ independent of $T$. It follows
that there is a positive constant $c_4$ independent of $T$ such
that
\[\frac{\s_k}{\s_l}(g)\ge c_4,\]
at point $(x_0,t_0)$. Then the Proposition follows, as $|u|$ is
bounded by Proposition \ref{pro2}. \qed

\medskip

\noindent{\it Proof of Theorem \ref{mainthm}.} Now we can prove
Theorem \ref{mainthm}. First, from Propositions \ref{pro2} and
\ref{pro3} we have $T^*=\infty$. Then, by Krylov's Theorem
\cite{Krylov}, the flow has $C^{2,\a}$ estimates. (\ref{4})
implies that
\[\int_0^\infty\int_M ({\s_k(g)}-r_{k,l}\s_l(g))^2 dg dt
<\infty,\]
which, in turn, implies that there is a sequence
$\{t_l\}$ such that
$$\int_M ( {\s_k(g)}-r_{k,l}\s_l(g))^2(t_l) dg \to 0$$ as
$t_l \to \infty$. The above estimates imply
that $g(t_l)$ converges in $C^{2,\a}$ to a conformal
metric $h$, which is a solution of (\ref{quot-eq}).

Now we want to use Simon's argument \cite{Simon} to prove that $h$
is the unique limit of flow (\ref{flow}) as in \cite{GW2} (see
also \cite{Andrews}). Since the arguments are essentially the
same, here we only give a sketch. First, with the $C^{2,\alpha}$
regularity estimates (and higher regularity estimates follows from
the standard parabolic theory) established for flow (\ref{flow}),
one can show that,
\[\lim_{t \to \infty} \|\frac{\s_k}{\s_l}(g(t))-\beta\|_{C^{4,\a}(M)}=0, \]
for some positive constant $\beta$. It is clear that
$\frac{\s_k}{\s_l}(h)=\beta$. By Proposition \ref{pro3} and the
Newton-MacLaurin inequality, there is a constant $c>1$ such that
$c^{-1}\le \s_l(g(t))\le c$. We want to show that flow
(\ref{flow}) is a pseudo-gradient flow, though it is not a
gradient flow. The crucial step is to establish the    ``angle
estimate" (\ref{gradient}) for the $L^2$ gradient of some proper
functionals. The estimate (\ref{gradient}) enables one to conclude
that the flow will converge to a unique stationary solution, we
refer to \cite{GW2} for the detailed argument.

Now we may switch the background metric to $h$ and all derivatives
and norms are taken with respect to the metric $h$ (since we have
all the a priori estimates for $h$).

Here we only give a proof for $l<k<n/2$. The proof for the other cases
is similar after taking consideration of the corresponding functionals.

Consider a functional defined by
\[{\FF}_{k,l}(g)=
 \left(\int \s_l(g)dg\right)^{-\frac {n-2k}{n-2l}}\int_M \s_k(g) dg.\]
Its $L^2$-gradient is
\[\n \F_{k,l}=-c_0((\s_k(g)-\tilde r_{k,l}(g) \s_l(g))e^{-nu},\]
where 
$c_0$ is a non-zero
constant and $\tilde r_{k,l}(g)$ is given by
 \[\tilde r_{k,l}(g):=\frac{\int_M \s_k(g) dg}
{\int_M \s_l(g) dg} ,\] which is different from $r_{k,l}$. But it
is easy to check that $r_{k,l}(t)-\tilde r_{k,l}(t)\to 0$ as $ t
\to \infty$. Since $\frac{\s_k}{\s_l}(g(t))$ is very close to a
constant for large $t$, from (\ref{4}) we have
\begin{equation} \label{gradient}\begin{array}{rcl}
\frac d {dt}\FF_{k,l}(g) &\le & \ds\vs -c \int_M
(\frac{\s_k}{\s_l}(g)-r_{k,l})
\left(\log \frac{\s_k}{\s_l}(g)-\log r_{k,l}\right)\s_l(g)dg\\
 &\le& \ds\vs -c \left(
 \int_M \left|\frac{\s_k}{\s_l}(g)-r_{k,l}\right|^2\s_l(g) dg
\int_M \left|\log \frac{\s_k}{\s_l}(g)-\log r_{k,l}\right|^2\s_l(g) dg\right)^{1/2}\\
  &\le& \ds\vs -c \left(
 \int_M |\frac{\s_k}{\s_l}(g)-r_{k,l}|^2\s_l(g) dg\right)^{1/2}
 \left(\int_M |\frac{dg}{dt}|^2 \s_l(g)dg\right)^{1/2} \\
   &\le& \ds\vs -c\left(
 \int_M |\frac{\s_k}{\s_l}(g)-\tilde r_{k,l}|^2\s_l(g) dg\right)^{1/2}
  \left(\int_M |\frac{dg}{dt}|^2\s_l(g)dg\right)^{1/2} \\
   &\le& \ds\vs -c\left(
 \int_M |{\s_k(g)}-\tilde r_{k,l}{\s_l(g)}|^2\s_l(g) dg\right)^{1/2}
  \left(\int_M |\frac{dg}{dt}|^2\s_l(g) dg\right)^{1/2} \\
  &\le& \ds\vs -c\left(
 \int_M |\n \FF_{k,l}|^2 dh\right)^{1/2}
  \left(\int_M |\frac{dg}{dt}|^2dh\right)^{1/2},
  \end{array}\end{equation}
where $c>0$ is a constant varying from line to line. The angle estimate
(\ref{gradient})
means that flow (\ref{flow}) is a pseudo-gradient flow. Finally from
(\ref{gradient}) we can
apply  the argument given in \cite{Simon} to show the uniqueness of the limit
as in \cite{GW2}. \qed

\section{Proof of Theorem \ref{inqs}}
{{From}} the flow approach we developed, we have
\begin{pro}\label{thm-all} Let $(M,g_0)$ be a compact, connected
and oriented \lcf manifold with $g_0\in \Gamma_k^+$ and $0 \le l<
k\le n$. Let $\tFF_k$ defined as in (\ref{tff}), then there is
$g_E \in {\cal C}_k$ satisfying equation (\ref{quot-eq}) such that
\[ 0<\tFF_k(g_E) \le \tFF_k(g),\]
for any $g\in {\cal C}_k$ with $\tFF_l(g_E) = \tFF_l(g).$
Moreover, if $(M,g_0)$ is conformally equivalent to a space form, then $(M,g_E)$
is also a space form.
\end{pro}

\pr  The case $l=0$ has been treated in Proposition \ref{GW}. We
may assume $l\ge 1$ in the rest of proof. When $(M,g_0)$ is
conformally equivalent to a space form, then any solutions of
(\ref{quot-eq}) are  metrics of constant sectional curvature (this
can be proved using the method of moving plane of \cite{GNN} as in
\cite{Jeff2}, see Corollary 1.1 in \cite{Yanyan} for more general
statement), and hence have the same $\tFF_k$ if they have been the
same $\tFF_l$. Hence the Proposition follows from Theorem
\ref{mainthm}.

Now we remain to consider the case $k<n/2$ and
$(M,g_0)$ is not conformally equivalent to a space form.
We will follow the same argument in the proof of Proposition \ref{GW}.
Here we need the local estimates in Theorem \ref{thm-local} for the quotient equation (\ref{quot-eq}).

First we want to show
\begin{equation}\label{beta0}
\inf_{g\in {\cal C}_k, \tFF_l(g)=1} \tFF_k(g) =: \beta_0>0.
\end{equation}
Suppose $\beta_0=0$.
By the result for flow (\ref{flow}), there is  a sequence
$g_i=e^{-2u_i}g_0\in {\cal C}_k$ with $\tFF_l(g_i)=1$
and
\[\frac{\sk}{\s_l}(g_i)={\b_i},\quad \lim_{i\to \infty}{\b_i}=0,\]
The scaled metric $\tilde g_i= e^{-2\tilde u_i}g_0$
with $\tilde {u}_i={u_i} -\frac 1{2(k-l)} \log {\b_i}$ satisfies
\begin{equation}\label{appl_12}
\frac{\s_k}{\s_l}(\uuult)= e^{-2(k-l)\tilde u_i}.\end{equation}
By Proposition 1,
\[Cvol(\tilde g_i)^{\frac{n-2l}n} \le \tFF_l(\tilde g_i)= {\b_i}^{\frac {n-2l}{2(k-l)}}
\to 0 \quad \hbox{ as } i\to \infty,\] Since $vol(\tilde g_i)$ is
tending to $0$, by Corollary \ref{newenergy},
 \[{m_i}:=\inf_M \tilde {u}_i \to +\infty \quad \text{as} \quad i\to \infty.\]

Now at the minimum point $x_i$ of $\tilde u_i$, by equation (\ref{appl_12}),
\[\frac{\s_k}{\s_l}(S_{g_0}) \le \frac{\s_k}{\s_l}(\uuult)= e^{-2(k-l)m_i}
\rightarrow 0.\] This is a contradiction to the fact  that $g_0
\in \Gamma_k^+$, proving (\ref{beta0}).

Finally we prove the existence of an extremal metric in this case.
{From} above argument, there is a minimization sequence $g_i \in {\cal C}_k$,
with $\tFF_l(g)=1$, and $\frac{\s_k(g_i)}{\s_l(g_i)}=\b_i,$
with $\b_i$ decreasing and bound below by a positive constant.
As $(M,g_0)$ is not conformally equivalent to ${\S^n}$ by assumption,
it follows from Theorem 1.3 in \cite{GLW} (see also \cite{LL3}) that
the metrics converge (by taking a subsequence) to some $g_E$ which attains
the infimum $C_S$.
\qed

\medskip

\noindent{\it Proof of  $(${\bf B}$)$ of Theorem \ref{inqs}.} The cases $l=n/2$
and $k=n/2$
were considered in \cite{GW2} and \cite{GVW}. Hence we assume that
$k\neq n/2$ and $l\neq n/2$.
 Let us consider
\[{\FF}_{k,l}(g)=
 \left(\int \s_l(g)dg\right)^{-\frac {n-2k}{n-2l}}\int_M \s_k(g) dg.\]
Since $\F_{k,l}$ is invariant under the transformation
$g$ to $e^{-2a}g$ for any constant $a$, Proposition \ref{thm-all} implies that
for any $g\in {\cal C}_k$
\[\FF_{k,l}(g) \le \FF_{k,l}(g_E)=:C(n,k,l).\]
It is clear that $C(n,k,l)$ depends only on $n,k,l$.

Let $c_0=\F_{k,l}(g_E)$.
{{From}} Proposition \ref{thm-all}, we have
\begin{equation}\label{f2}\begin{array}{rcl}
\int_M \s_k(g) dg &\le &C(n,k,l) \left(\int_M\s_l(g) dg \right)^
{\frac{n-2k}{n-2l}}\\
&=& C(n,k,l) \left(\int_M\s_l(g) dg \right)^{\gamma k}
\left(\int_M\s_l(g) dg \right)^{\frac k l},
\end{array}\end{equation}
where $\gamma=\frac {n-2k}{k(n-2l)}-\frac 1l.$
It is clear that $\gamma>0$ when $l>n/2$ and
$\gamma<0$ when $l<n/2$.

We first consider the case $l>n/2$. In this case, by Proposition \ref{GW} we have
\[\int_M\s_l(g) dg \le c_1 vol (g)^{\frac {n-2l} n},\]
where $c_1=\F_l(g_e)$. It follows that
\begin{equation}\label{f1} \left(\int_M\s_l(g) dg\right)^\gamma
\le c_0^\gamma vol (g)^{\frac{l-k}{kl}}.\end{equation}
Hence
\[\begin{array}{rcl}
\ds\vs (\F_k(g))^{1/k}& =& \ds \left(vol(g)^{-\frac{n-2k} n} \int_M \s_k(g) dg
\right)^{\frac 1k}\\
&\le &c_0^{\frac 1k} \ds\vs \left(vol(g)^{-\frac{n-2l} n} \int_M \s_l(g) dg
\right)^{\frac 1l}\\
&=& \ds c_0^{\frac 1k} (\F_l(g))^{1/l}.\end{array}\]
The equality holds if and only if $g$ is a metric of constant sectional curvature.

Consider the case $l<n/2$. In this case, by Proposition \ref{GW} again we have
\[\int_M\s_l(g) dg \ge c_1 vol (g)^{\frac {n-2l} n},\]
where $c_1=\F_l(g_e)$. Since $\gamma<0$, we have
(\ref{f1}). The same argument given in the previous case gives the same
conclusion.

Finally, since $k\ge n/2$, $(M,g_0)$ is conformally equivalent to
a space form (\cite{GVW}). The existence of the extremal metric
which attains the equality case follows from Proposition
\ref{thm-all}. And the constant $C(n,k,l)$ is easy to calculate.
\qed

\medskip

\noindent{\it Proof of $(${\bf C}$)$ of Theorem \ref{inqs}.}
Let us first consider the case $l<n/2$. Let $g\in {\cal C}_{n/2}$.
Choose $a$ such that
$\int_M\s_l(e^{-2a }g)dvol (e^{-2a} g)=\int_M\s_l(g_0)dg_0.$
It is easy to see that
\[a=\frac 1 {n-2l} \{\log \int_M\s_l(g)dg-\log\int_M\s_l(g_0)dg_0\}.\]
By   Proposition \ref{thm-all}, we have
\begin{eqnarray*}
\En(g) &=& \En(e^{-2a}g)+a\int_M \s_{n/2} (g) dg\\
& \ge & a \int_M \s_{n/2} (g_0) dg_0 \\
&=& \frac 1 {n-2} \int_M \s_{n/2} (g_0) dg_0\left\{
\log \int_M\s_l(g)dg-
\log\int_M\s_l(g_0)dg_0\right\}.\end{eqnarray*}
This proves the Theorem
for the case $l<n/2$.

Now we consider the case $l>n/2$.
 \ref{thm-all}. For any $g\in {\cal C}_l$
we choose
\[a=(\int_M \s_{n/2}(g) dg)^{-1} \En(g)\]
 such that $\En(e^{-2a}g)=\En(g_0)$.
 Recall that $\tFF_{n/2}=\En$. By Proposition \ref{thm-all} again, we have
 \begin{eqnarray*}\tFF_l(g) &=&
 \frac1 {n-2l } \int_M\s_l(g) dg\\
 &=& \frac1 {n-2l}e^{-(2l-n)a}\int_M\s_l(e^{-2a}g)
 dvol(e^{-2a}g)\\
 &\ge& \frac1 {n-2l}e^{-(2l-n)a}\int_M\s_l(g_0)
 dg_0\\
 &=& \frac1 {n-2l} \exp\left\{(n-2l)\left(\int_M \s_{n/2}(g) dg\right)^{-1} \En(g)\right\}
 \int_M\s_l(g_0)
 dg_0.\end{eqnarray*}
 Since $(M,g_0)$ is conformally equivalent to a space form in this case, the existence
of the extremal metric can be proved along the same line as in
part ({\bf B}) of the Theorem. Note that since $n$ is even,
$(M,g_0)$ is the standard sphere. The computation of $C_{MT}$ is
straightforward.\qed

\medskip

\noindent{\it Proof of $(${\bf A}$)$ of Theorem \ref{inqs}.}
Inequality (\ref{Sob-ineq}) follows from (\ref{beta0}) in the
proof of Proposition \ref{thm-all}. The existence of the extremal
metric has also proved there. As for inequality (\ref{best1}),
since its proof is of different spirit and inequality itself is of
independent interest, we will devote it in the next section
(Theorem \ref{best}). The constant $C_{S,k,l}({\S^n})$ in
(\ref{best1}) can be computed easily. \qed

\section{The best constant}

In this section, we address the question of the best constant in part
({\bf A})
of Theorem \ref{inqs}.
As in the Yamabe problem (i.e., $k=1$ and $l=0$), for $0\le l< k <n/2$ we define
\[Y_{k,l}(M,[g_0])=\inf_{g\in {\cal C}_k} (\F_{l}(g))^{-\frac {n-2k}{n-2l}}
\F_k(g)=\inf_{g\in {\cal C}_k} (\int_M\s_l(g) dg)^{-\frac {n-2k}{n-2l}}
\int_M\s_k(g) dg.\]
It is clear that $Y_{k,l}(M,[g_0])=C_s^{n-2k}$. In this section we prove
\begin{thm}\label{best} For any compact, oriented \lcf manifold $(M,g_0)$,
we have
\begin{equation}\label{best11} Y_{k,l}(M,[g_0]) \le Y_{k,l}(\S^n, g_{\S^n}),\end{equation}
where $g_{\S^n}$ is the standard metric of the unit sphere.
\end{thm}
When $k=1$ and $l=0$, this
was proven by Aubin (e.g., see \cite{Aubin}) for general compact
manifolds.
To prove Theorem \ref{best} we need to construct a
sequence of ``blow-up" functions which belong to ${\cal C}_k$. This is a
delicate part of the problem.

We need two Lemmas.
 \begin{lem}\label{lemG1}
 Let $D$  be the unit disk in $\R^n$ and $ds^2$ the standard Euclidean metric.
 Let $g_0=e^{-2u_0}ds^2$
be a metric on $D$ and $g_0\in\Gamma_k$ with $k<n/2$. Then there
is a conformal metric $g=e^{-2u}ds^2$ on $D\backslash\{0\}$ of
 positive $\Gamma_k$-curvature with the following properties:
 \begin{itemize}
 \item [1).] $\s_k(g)>0$ in $D\backslash\{0\}$.
 \item [2).] $u(x)=u_0(x)$ for $r=|x| \in (r_0,1]$.
 \item [3).] $u(x)=a+\log r$
 for $r=|x|\in (0, r_3)$ and some constant $a$.
 \end{itemize}
 for some constants $r_0$ and $r_3$ with $0<r_3<r_0<1$.
 \end{lem} \pr
 Let $v$ be a function on $D$ and $\tilde g=e^{-2v}g_0$.
 By the transformation formula of
 the Schouten tensor, we have
\begin{equation}\label{g1}\begin{array}{rcl}
S(\tilde g)_{ij}&=& \ds\vs \n_{ij}^2(v+u_0) +\n_i(v+u_0)\n_j(v+u_0) -\frac 12
|\n(v+u_0)|^2 \d_{ij} \\
&=& \ds\vs \n^2_{ij} v+\n_i v\n_j v+\n_i v \n_j u_0 +\n_j v\n_i u_0\\
&&  \ds +(\frac 12|\n v|^2+ \n v \n u_0)\d_{ij}+S_{g_0} \end{array}
\end{equation}
Here $\n $ and $\n^2$ are the first and the second derivatives with respect
to the standard metric $ds^2$.
Let $r=|x|$. We want to find a  function $v=v(r)$
with $\tilde g \in \Gamma_k^+$  and
\[ v' =\frac{\a(r)} r,\]
where $\a=1$ near $0$ and $\a=0$ near $1$. From (\ref{g1})
we have
\begin{eqnarray} \label{g1.1}
S({\tilde g})_{ij}&= &
\frac {2\a -\a^2}{2r^2} \d_{ij}+\left( \frac {\a'}{r} +\frac{\a^2-2\a} {r^2}
\right) \frac
{x_ix_j}{r^2}+ S({g_0})_{ij}+ O(|\n u_0|) \frac{\a}r ,\end{eqnarray}
where $O(|\n u_0|)$ is a term bounded by a constant $C_1$ depending
only on $ \max |\n u_0|$. Let $A(r)$ be an $n\times n$ matrix with entry
$a_{ij}=S(\tilde g)_{ij}-S(g_0)_{ij}$. Hence
\[\s_k(\tilde g)=e^{-2k(v+u_0)}
\s_k\left(A+ S({g_0})\right).\] To our aim, we need to
find $\a$ such that $A+S(g_0) \in \Gamma_k^+$.
Let $\e\in (0,1/2)$ and $r_0= \min\{\frac 12, C_1 \e\}$.
We will choose $\a$ such that
\begin{equation}
\label{eq-g1}\a (r) \in [0,1] \text{ and } \a(r)=0, \text{ for } r\in [r_0, 1].
\end{equation}
Since $\s_k(\tilde g)=e^{2k(v+u_0)}\s_k(A(r)+S(g_0))$,
we want to find $\a$ such that $\s_k(A(r)+S(g_0))>0$.
It is clear to see that for $r\in [0,r_0]$
\[A(r) \ge \left (\frac {2\a -\a^2-\e \a}{2r^2} \d_{ij}
+\left( \frac {\a'}{r} +\frac{\a^2-2\a} {r^2}
\right) \frac{x_ix_j}{r^2} \right),\]
as a matrix. This implies that
\begin{equation}\label{5.1}\s_k(A(r))\ge \frac{(n-1)!}{k!(n-k)!}
\left(\frac {2\a-\a^2-\e\a}{2 r^2}\right)^k
\left(n-2k+2\frac {r\a'-\e\a}{2\a-\a^2-\e\a}\right).\end{equation}
One can easily check that for any small $\d>0$,
\begin{equation}\label{H1}\a(r)=\frac{2(1-\e)\d}{\d+r^{\frac{1-\e}2}}
\end{equation}
is a solution of
\[(2-\e)\a-\a^2=-4(r\a'-\e\a).\]

Now we can finish our construction of $\a$.
Since $S(g_0)\in\Gamma_k^+$, by the openness of $\Gamma_k^+$ we can
choose $r_1\in (0,r_0)$ and an non-increasing function
$\a:[r_1,r_0] \subset [0, 1)$ such that $\s_k(\tilde g)>0$ and
$\a(r_1)>0$. Now  we choose
a suitable $\d>0$  and $\a$ in the form (\ref{H1}).Then find
  $r_2\in (0,r_1)$ with $\a(r_2)=1$.
 It is clear that $\s_k(A(r))>0$ on
$[r_2,r_1]$.
Define $\a(r)=1$ on $[0,r_2]$. We may smooth $\a$
such that the new resulted conformal metric $g$ satisfying all
conditions in Lemma \ref{lemG1}. \qed

\begin{rem} From Lemma \ref{lemG1}, one can prove that the connected sum of
 two \lcf manifolds  $(M_1,g_1)$ and $(M_2,g_2)$ with $g_1, g_2\in
 \Gamma_k^+$ ($k<n/2)$ admits a metric in $\Gamma_k^+$.
This is also true for general manifolds. Namely, the connected sum of
 two compact manifolds  $(M_1,g_1)$ and $(M_2,g_2)$ with $g_1, g_2\in
 \Gamma_k^+$ ($k<n/2)$ admits a metric in $\Gamma_k^+$. The proof is given in
 \cite{GLW2}.\end{rem}
\begin{lem}\label{lemG2}For
any small constants ${\d}>0$ and $\e>0$, there exists
a function $u:\R^n\backslash\{0\}\to 0$ satisfying:
\begin{itemize}
\item[1.] The metric $g=e^{-2u}dx^2$ has positive $\Gamma_k$-curvature.
\item[2.] $u=\log(1+|x|^2)+b_0$ for $|x|\ge {\d}$, i.e.,
$(\{x\in \R^n\,|\, |x|\ge {\d}\}, g)$ is a part of a sphere.
\item[3.] $u=\log|x|$  for $|x|\le \d_1$, i.e.,
$(\{x\in \R^n\,|\, 0<|x|\le \d_1\}, g)$ is a cylinder.
\item[4.] $vol(B_{{\d}}\backslash B_{\d_1}, g) \le C \d^{-\frac{2n}{1-\e_0}}.$
\item[5.]
$\int_{B_{{\d}}\backslash B_{\d_1}}\s_k(g) dvol(g)\le
 C \d^{-\frac{2(n-2k)}{1-\e_0}}$, for any $k<n/2$,\end{itemize}
where $C$ is a constant independent of ${\d}$,
$\d_1=\d^{\frac{3-\e_0}{1-\e_0}}$ and
 $b_0\sim\frac{3-\e_0}{1-\e_0} \log\d $.\end{lem}

\pr Let ${\d}\in (0,1)$ be any small constant. For any small constant
$\e_0>0$, we define $u$ by
\[
u(r)=\begin{cases}
\log(1+r^2)+b_0, & r\ge {\d}\\
 \ds\vs -\frac 2{1-\e_0}
\log \frac {1+\d^{3-\e_0} r^{-(1-\e_0)}}{2}
+\frac {3-\e_0}{1-\e_0}\log\d
& r\in ({\d_1}, \d)\\
\ds\vs \log  r, &
r\le \d_1,\end{cases},\]
where  $\d_1=\d^{\frac{3-\e_0}{1-\e_0}}$ and
\[b_0=-\log (1+{\d}^2)
-\frac 2{1-\e_0}
\log \frac {1+\d^2} 2
+\frac {3-\e_0}{1-\e_0}\log\d.\]
As in the proof of Lemma 5, we write $u'(r)=\frac {\a(r)} r$. It is
easy to see that
$\a:\R_+\to \R_+$ by
\[
\a(r)=\begin{cases} \ds\vs
\frac{2r^2}{1+r^2}, & r\ge {\d},\\
\ds\vs  \frac{2\d^{3-\e_0}}
{{\d^{3-\e_0}}+ r^{1-\e_0}}, & r\in (\d_1, \d),\\
1, & r\le
\d_1.\end{cases}\]
One can check all conditions in the Lemma, except the smoothness of $u$,
which is $C^{1,1}$.
We first check (1). By a direct computation, see for example (13),
we have
\[\s_k(e^{-2u}|dx|^2)= e^{2k u(r)}
\frac{(n-1)!}{k!(n-k)!}
\left(\frac {2\a-\a^2}{2 r^2}\right)^k
\left(n-2k+2\frac {r\a'}{2\a-\a^2}\right).\]
In the interval $({\d_1}, {\d})$, $\a \in (0, 2)$ satisfies
\[\frac {2r\a'}{2\a-\a^2}=-(1-\e_0).\]
Since $k<n/2$, we have $\s_k(e^{-2u}|dx|^2)>0$.
One can also directly to check (4) and (5). Here we only check (5).
A direct computation gives
\begin{eqnarray*}
\ds\vs \int_{B_{{\d}}\backslash B_{\d_1}} \s_k(g) dvol(g) & \le &
c\ds \int_{\d_1}^{\d} e^{-(n-2k)u(r)} r^{-2k} r^{n-1}dr\\
&\le & c{\d}^{-(n-2k)\frac{3-\e_0}{1-\e_0}}
\int_{\d_1}^{{\d}} r^{n-2k-1}dr \\
& \le &c \ds\d^{-\frac{2(n-2k)}{1-\e_0}}.\end{eqnarray*}
{From} our construction, we only have $u\in C^{1,1}$. But, for
$\delta>0$ fixed, we can smooth $\a$ so that $u\in C^{\infty}$
satisfies all conditions (1)-(5). \qed

\noindent{\it Proof of Theorem \ref{best}.}
Let $p\in M$ and $U$  a neighborhood of
$p$ such that $(U, g)$ is conformally flat, namely $(U,g)=(D, e^{-2u_0}|dx|^2)$.
Applying Lemma \ref{lemG1}, we obtain a conformal metric $u$ satisfying
conditions 1)-3) in Lemma \ref{lemG1} with constants ${r_0}, r_3$ and $a$.
By adding a constant we may assume $a=0$. Now applying Lemma \ref{lemG2}
 for any small constant ${\d}>0$ we
have a conformal metric $g_\d=e^{-2u_\d}|dx|^2$ on $\R^n\backslash\{0\}$.
Consider the rescaled function
\[\tilde u_\d=u_\d(\frac {\d_1}{r_3} x)-\log \frac {\d_1}{r_3}.\]
Now $u$ and $\tilde u_\d$ are the same in $\{0<|x|<r_3\}$.
Consider the following conformal transformation
\[f(x)= \frac{r^2_3} 2\frac{x}{|x|^2},\]
which maps $\{r_3/2\le |x|\le r_3\}$ into itself and
maps one of boundary components to another with opposite orientations.
Now we define a new function on $M$
by
\[w_\d(x)=\begin{cases}
\vs 0, & |x| \ge r_0 ,\\
\vs u-u_0 ,& r_3/2\le |x| \le r_0 ,\\
 \tilde u_\d(f(x))+2\log|x| -\log \frac {r^2_3}2-u_0, &|x| \le r_3/2. \end{cases}\]
Since  $u$ and $\tilde u_\d$ are the same in $\{0<|x|<r_3\}$, it
clear that $w_\d(x)$ is smooth on $M$. Consider the conformal
metric $g_\d=e^{-2w_\d}g$ and compute, using Lemmas \ref{lemG1}
and \ref{lemG2}
\begin{eqnarray*}
 \int_M \s_k(g_\d)dvol(g_\d)& = &\int_{\{|x|\le r_3/2\}}\s_{k}(g_\d) d{g_\d}+O(1)
 \\
 &=&e^{-(n-2k)b_0}\int_{\R^n\backslash\{|x|\le \d\}}
\s_k(g_{\S^n})dvol(g_{\S^n})+O(1)\d^{-\frac{2(n-2k)}{1-\e_0}}\\
&=&\d^{-\frac{3-\e_0}{1-\e_0}(n-2k)}\int_{\R^n\backslash\{|x|\le \d\}}
\s_k(g_{\S^n})dvol(g_{\S^n})+o(\d^{-\frac{3-\e_0}{1-\e_0}(n-2k)})\end{eqnarray*}
and
\[\int_M \s_l(g_\d)dvol(g_\d)=\d^{-\frac{3-\e_0}{1-\e_0}(n-2l)}vol(\R^n\backslash\{|x|\le \d\}, g_{\S^n})
+o(\d^{-\frac{3-\e_0}{1-\e_0}(n-2l)}),\]
where $g_{\S^n}=\frac 1{(1+|x|^2)^2}|dx|^2$ is the standard metric
of the sphere and $O(1)$ is a term bounded by a constant independent of
$\d$. Now it is readily to see
\[Y_{k,l}(M) \le \lim_{{\d}\to 0} Y_{k,l}(g_{{\d}}) \to Y_{k,l}(\S^n),\]
as ${\d}\to 0$. \qed

\end{document}